\documentclass[leqno,11pt]{article}  
\usepackage{amsmath}
\usepackage{amssymb}

\usepackage[german,english]{babel}
\usepackage{amsthm}

\title{Finiteness of de Rham cohomology in rigid analysis}

\author{\textsc{Elmar Grosse-Kl\"onne}}
\date{}

\theoremstyle{plain} 
\newtheorem{satz}{Theorem}[section]  
\newtheorem{kor}[satz]{Corollary}  
\newtheorem{lem}[satz]{Lemma}  
\newtheorem{pro}[satz]{Proposition}  
\newcommand{\ho}{\mbox{\rm Hom}}  
\newcommand{\ext}{\mbox{\rm Ext}} 
\newcommand{\spec}{\mbox{\rm Spec}}  

\newcommand{\spm}{\mbox{\rm Sp}}  
\newcommand{\spf}{\mbox{\rm Spf}}  

\newcommand{\codim}{\mbox{\rm codim}}  
\newcommand{\bi}{\mbox{\rm im}}  
\newcommand{\ke}{\mbox{\rm Ker}}  

\newcommand{\dr}{\mbox{\rm DR}}  
\newcommand{\q}{\mbox{\rm Frac}}  
\newcommand{\der}{\mbox{\rm D}}

\newcommand{\kodi}{\mbox{\rm coh.dim}}
\newcommand{\sym}{\mbox{\rm Sym}}

\theoremstyle{remark}

\theoremstyle{definition}

\begin{document}
\maketitle
\footnote[0]{I wish to thank Peter Schneider who patiently read early drafts of my attempts on this subject and gave helpful comments. I am also very  grateful to Pierre Berthelot for an invitation to the University of Rennes I, with discussions of decisive importance during this visit, and for further remarks on the text later on. I thank Annette Huber-Klawitter and Klaus K\"unnemann for answering me various questions. Finally, I thank the referee for his careful reading and helpful advice.}

\begin{abstract} 

For a big class of smooth dagger spaces --- dagger spaces are 'rigid spaces with overconvergent structure sheaf' --- we prove finite dimensionality of de Rham cohomology. This is enough to obtain finiteness of Berthelot's rigid cohomology also in the non-smooth case. We need a careful study of de Rham cohomology in situations of semi-stable reduction.

\end{abstract}

%


\begin{center} {\bf Introduction} \end{center}

Let $(R,\pi,k=\q(R),\bar{k}=R/(\pi))$ be a complete discrete valuation ring of mixed characteristic. It is a simple observation that the de Rham cohomology $H_{dR}^*(X)$ of a positive dimensional smooth affinoid $k$-rigid space $X$ computed with respect to its (usual) structure sheaf is not finite dimensional. The idea of instead using an overconvergent structure sheaf arises naturally from the paper \cite{mw} of Monsky and Washnitzer: 
The Monsky-Washnitzer cohomology of a smooth affine $\bar{k}$-scheme $\spec(A)$ is the de Rham cohomology of $\tilde{A}^{\dagger}\otimes_Rk$, where $\tilde{A}^{\dagger}$ is a weakly complete formal lift of $A$. Monsky-Washnitzer cohomology has recently been shown to be finite dimensional (independently by Berthelot \cite{berfi} and Mebkhout \cite{mebend}). The algebra $\tilde{A}^{\dagger}\otimes_Rk$ can be geometrically interpreted as a $k$-algebra of overconvergent functions on the rigid space $\spm(\tilde{A}\otimes_Rk)$, where $\tilde{A}$ is a lifting of $A$ to a formally smooth $\pi$-adically complete $R$-algebra. In \cite{en1dag} we introduce a category of '$k$-rigid spaces with overconvergent structure sheaf', which we call $k$-dagger spaces, and study a functor $X\mapsto X'$ from this category to the category of $k$-rigid spaces which is not far from being an equivalence. 
For example, $X$ and $X'$ have the same underlying $G$-topological space and the same stalks of structure sheaves. Finiteness of Monsky-Washnitzer cohomology implies finiteness of de Rham cohomology for affinoid $k$-dagger spaces with good reduction --- in the above notation, the algebra $\tilde{A}^{\dagger}$ gives rise to the affinoid $k$-dagger space $X$ with $\Gamma(X,{\cal O}_X)=\tilde{A}^{\dagger}\otimes_Rk$. Our main result generalizes this as follows:\\

{\bf Theorem A:} (=Corollary \ref{ende} $+$ Theorem \ref{kokomend}) {\it Let $X$ be a quasi-compact smooth $k$-dagger space, $U\subset X$ a quasi-compact open subset, $Z\to X$ a closed immersion. Then $T=X-(U\cup Z)$ has finite dimensional de Rham cohomology $H_{dR}^*(T)$.}\\

By \cite{en1dag}, 3.2, this implies finiteness of de Rham cohomology also for certain smooth $k$-rigid spaces $Y$: For example, if $Y$ admits a closed immersion $i$ into a polydisk without boundary (at least if $i$ extends to a closed immersion with bigger radius), or if $Y$ is the complement of a quasi-compact open subspace in a smooth proper $k$-rigid space. But our main corollary is of course:\\

{\bf Corollary B:} (=Corollary \ref{rifi}) {\it For a ${\bar k}$-scheme $X$ of finite type the $k$-vector spaces $H_{rig}^*(X/k)$ (see \cite{berfi}) are finite dimensional.}\\

We do not reprove finiteness of Monsky-Washnitzer cohomology, rather we reduce our Theorem A to it. A big part is devoted to the study of de Rham cohomology in situations of semi-stable reduction. We need and prove:\\

{\bf Theorem C:} (=Theorem \ref{monovor}) {\it Let ${\cal X}$ be a strictly semi-stable formal $R$-scheme, let ${\cal X}_{\bar{k}}=\cup_{i\in I}Y_i$ be the decomposition of the closed fibre into irreducible components. For $K\subset I$ set $Y_K=\cap_{i\in K}Y_i$. Let $X^{\dagger}$ be a $k$-dagger space such that its associated rigid space is identified with ${\cal X}_k$. For a subscheme $Y\subset{\cal X}_{\bar{k}}$ let $]Y[_{\cal X}^{\dagger}$ be the open dagger subspace of $X^{\dagger}$ corresponding to the open rigid subspace $]Y[_{\cal X}$ of ${\cal X}_k$. Then for any $\emptyset\ne J\subset I$ the canonical map $H_{dR}^*(]Y_J[_{\cal X}^{\dagger})\to H_{dR}^*(]Y_J-(Y_J\cap(\cup_{i\in I-J}Y_i))[^{\dagger}_{\cal X})$ is bijective.}\\ 

Another important tool is de Jong's theorem on alterations by strictly semi-stable pairs, in its strongest sense.\\
We proceed as follows. After recalling some facts on dagger spaces in section 0, we formulate in section 1 some basic concepts about ${\cal D}$-modules on rigid and dagger spaces. This follows the complex analytic case, see for example \cite{mebdmod}. Instead of reproducing well known arguments, we focus only on what is specific to the non-archimedean case. Then we construct a long exact sequence for de Rham cohomology with supports in blowing-up situations. As in \cite{hadr} (for algebraic $k$-schemes), it results from the existence of certain trace maps for proper morphisms; we define such trace maps, based on constructions from \cite{chia}, \cite{bey}, \cite{vdpse}. Finally we prove the important technical fact that the de Rham cohomology $H_{dR}^*(X)$ of a smooth dagger space $X$ depends only on its associated rigid space $X'$; hence knowledge of $X'$ (for example a decomposition into a fibre product) gives information about $H_{dR}^*(X)$. In the section 2, we begin to look at $R$-models of the associated rigid spaces; more specifically, we consider the case of semi-stable reduction. Its main result is Theorem C. It enables us to reduce Theorem A, in the case where $U=\emptyset$ and $X$ has semi-stable reduction, to the finiteness of Monsky-Washnitzer cohomology. In section 3 we first prove Theorem A in case $U=\emptyset=Z$: After reduction to the case where $X$ is affinoid and defined by polynomials, we apply de Jong's theorem to an $R$-model of a projective compactification of $X$ to reduce to the finiteness result of section 2. The case of general $Z$ is handled by a resolution of singularities (\cite{bier}). Then we treat the case of general $U$ by another application of de Jong's theorem. The formal appearence of these last arguments bears similarity with the finiteness proofs in \cite{berfi}, \cite{hadr}. But there are also distinctive features: in the simultanuous control of special and generic fibre; and in particular in our second application of de Jong's theorem: We apply it to a certain closed immersion of $R$-schemes $\bar{X}_{\bar{k}}\cup\bar{Y}\to\bar{X}$, where the space $X-U$ we are interested in is realized in (the tube $]\bar{Y}[$ of) {\it the compactifying divisor, not} in its open complement.\\

\addtocounter{section}{-1}
\section{Dagger spaces}

Let $k$ be field of characteristic $0$, complete with respect to a non-archimedean valuation $|.|$, with algebraic closure $k_a$, and let $\Gamma^*=|k_a^*|=|k^*|\otimes\mathbb{Q}$.\\ 
We gather some facts from \cite{en1dag}. For $\rho\in\Gamma^*$ the $k$-affinoid algebra $T_n(\rho)$ consists of all series $\sum a_{\nu}X^{\nu}\in
k[[X_1,\ldots,X_n]]$ such that $|a_{\nu}|\rho^{|\nu|}$ tends to zero if $|\nu|\to\infty$. The algebra $W_n$ is defined to be $W_n=\cup_{\stackrel{\rho>1}{\rho\in\Gamma^*}}T_n(\rho)$ \footnote{The notation $W_n$ is taken from [{\it U. G\"unzer}, Modellringe in der nichtarchimedischen Funktionentheorie, Indag. math. {\bf 29} (1967), 334-342]. There the author assigns only the name of Washnitzer to this algebra. However, the referee pointed out that the name Monsky-Washnitzer algebra is the usual one.}. A $k$-dagger algebra $A$ is a quotient of some $W_n$; a surjection $W_n\to A$ endows it with a norm which is the quotient of the Gauss norm on $W_n$. All $k$-algebra morphisms between $k$-dagger algebras are continuous with respect to these norms, and the completion of a $k$-dagger algebra $A$ is a $k$-affinoid algebra $A'$ in the sense of \cite{bgr}. There is a tensor product $\otimes^{\dagger}_k$ in the category of $k$-dagger algebras.
As for $k$-affinoid algebras, one has for the set $\spm(A)$ of maximal ideals of $A$ the notions of rational and affinoid subdomains, and for these the analogue of Tate's acyclicity theorem (\cite{bgr},8.2.1) holds. The natural map $\spm(A')\to\spm(A)$ of sets is bijective, and via this map the affinoid subdomains of $\spm(A)$ form a basis for the strong $G$-topology on $\spm(A')$ from \cite{bgr}. 
Imposing this $G$-topology on $\spm(A)$ one gets a locally $G$-ringed space, an affinoid $k$-dagger space. (Global) $k$-dagger spaces are built from affinoid ones precisely as in \cite{bgr}. The fundamental concepts and properties from \cite{bgr} translate to $k$-dagger spaces.\\There is a faithful functor from the category of $k$-dagger spaces to the category of $k$-rigid spaces, assigning to a $k$-dagger space $X$ a $k$-rigid space $X'$ (to which we will refer as the associated rigid space; but we will use the notation $(?)'$ not only for this functor). $X$ and $X'$ have the same underlying $G$-topological space and the same stalks of structure sheaf. A smooth $k$-rigid space $Y$ admits an admissible open affinoid covering $Y=\cup V_i$ such that $V_i=U_i'$ for uniquely determined (up to non canonical isomorphisms) affinoid $k$-dagger spaces $U_i$. 
Furthermore, this functor induces an equivalence between the respective subcategories formed by partially proper spaces (see below). In particular, there is an analytification functor from $k$-schemes of finite type to $k$-dagger spaces. For a smooth partially proper $k$-dagger space $X$ with associated $k$-rigid space $X'$, the canonical map $H_{dR}^*(X)\to H_{dR}^*(X')$ between the de Rham cohomology groups is an isomorphism: this follows from applying \cite{en1dag}, 3.2 to the morphism between the respective Hodge-de Rham spectral sequences. By a dagger space not specified otherwise, we will mean a $k$-dagger space, and similarly for dagger algebras, rigid spaces etc..\\
In the sequel, all dagger spaces and rigid spaces are assumed to be quasi-separated. We denote by ${\bf D}=\{x\in k; |x|\le1\}$ (resp. ${\bf D}^0=\{x\in k; |x|<1\}$) the unitdisk with (resp. without) boundary, with its canonical structure of $k$-dagger or $k$-rigid space, depending on the context. For $\epsilon\in\Gamma^*$, the ring of global functions on the polydisk $\{x\in k^n;\mbox{ all }|x_i|\le\epsilon\}$, endowed with its canonical structure of $k$-dagger space, will be denoted by $k<\epsilon^{-1}.X_1,\ldots,\epsilon^{-1}.X_n>^{\dagger}$. The dimension $\dim(X)$ of a dagger space $X$ is the maximum of all $\dim({\cal O}_{X,x})$ for $x\in X$. We say $X$ is pure dimensional if $\dim(X)=\dim({\cal O}_{X,x})$ for all $x\in X$.\\
A morphism $f:X\to Y$ of rigid or dagger spaces is called  partially proper (cf. \cite[p.59]{huet}), if $f$ is separated, there is an admissible open affinoid covering $Y=\cup_i Y_i$, and for all $i$ admissible open affinoid coverings $f^{-1}(Y_i)=\cup_{j\in J_i} X_{ij}=\cup_{j\in J_i} X_{ij}'$ with $X_{ij}\subset\subset_{Y_i}X_{ij}'$ for every $j\in J_i$ (where $\subset\subset_{Y_i}$ is defined as in \cite{bgr}).\\

\begin{lem}\label{hironaka} Let $Z\to X$ be a closed immersion into an affinoid smooth dagger space. There is a proper surjective morphism $g:\tilde{X}\to X$ with $\tilde{X}$ smooth,  $g^{-1}(Z)$ a divisor with normal crossings on $\tilde{X}$, and $g^{-1}(X-Z)\to(X-Z)$ an isomorphism. 
\end{lem}

{\sc Proof:} Write $X=\spm(W_n/I), Z=\spm(W_n/J)$ with ideals $I\subset J\subset W_n$. Since these ideals are finitely generated, there is a $\rho>1$ and ideals $I_{\rho}\subset J_{\rho}\subset T_n(\rho)$ such that $I=I_{\rho}.W_n$ and $J=J_{\rho}.W_n$, and such that the rigid space $X_{\rho}=\spm(T_n(\rho)/I_{\rho})$ is smooth. Apply \cite[1.10]{bier} to the closed immersion $Z_{\rho}=\spm(T_n(\rho)/J_{\rho})\to X_{\rho}$ to get a morphism of rigid spaces $\tilde{X}_{\rho}\to X_{\rho}$ with the desired properties. Its restriction to the partially proper open subspace $\cup_{\rho'<\rho}\spm(T_n(\rho')/(I_{\rho}))\subset X_{\rho}$ is a morphism of partially proper spaces (compositions of partially proper morphisms are partially proper, \cite{huet}), hence by \cite{en1dag},2.27 is equivalent with a morphism of dagger spaces. The restriction of the latter to $X$ does the job.\\  

\section{${\cal D}$-modules}  

\addtocounter{satz}{1}{\bf \arabic{section}.\arabic{satz}}  For a smooth dagger (or rigid) space $X$ let $I=\ke({\cal O}_{X\times X}\to{\cal O}_X, f\otimes g\mapsto fg)$. 
For $n\in\mathbb{N}$ we view ${\cal P}_X^n={\cal O}_{X\times X}/I^{n+1}$ via $d_{1,n}:{\cal O}_X\to{\cal P}_X^n,\quad f\mapsto f\otimes1$ as ${\cal O}_X$-algebra and define ${\cal D}_X^n={\underline\ho}_{{\cal O}_X}({\cal P}_X^n,{\cal O}_X)$. 
Then ${\cal D}_X=\lim_{\stackrel{\to}{n}}{\cal D}_X^n$ comes in a natural way with the structure of a sheaf of rings on $X$ (as in \cite{EGA},IV,16.8.10). 
By modules over ${\cal D}_X$ which are not explicitly declared as ${\cal D}_X$-right-modules we always mean ${\cal D}_X$-left-modules. ${\cal O}_X$ becomes a ${\cal D}_X$-module setting $P.f=P(d_{2,n}(f))$ for $P\in{\underline\ho}_{{\cal O}_X}({\cal P}_X^n,{\cal O}_X), f\in{\cal O}_X$ and $d_{2,n}:{\cal O}_X\to{\cal P}_X^n, f\mapsto 1\otimes f$. 
For a morphism $X\to Y$ of smooth spaces with $J=\ke({\cal O}_{X\times_YX}\to{\cal O}_X, f\otimes g\mapsto fg)$ we define $\Omega_{X/Y}^1=J/J^2$. By means of $d:{\cal O}_X\to\Omega^1_{X/Y},\quad a\mapsto(1\otimes a-a\otimes1)$ we form in the usual way the complex $\Omega_{X/Y}^{\bullet}=(\bigwedge^{\bullet}\Omega_{X/Y}^1,d)$. 
If $Y=\spm(k)$ we write $\Omega_X^{\bullet}$. If $X$ is of pure dimension $r$ we write $\omega_X=\Omega_X^r$, and if $f:X\to Y$ is a morphism of smooth and pure dimensional spaces, we write $\omega_{X/Y}=\omega_X\otimes_{{\cal O}_X}f^*\omega_Y^{\otimes(-1)}$.\\

\addtocounter{satz}{1}{\bf \arabic{section}.\arabic{satz}}  Suppose $X=\spm(A)$ with a regular $k$-dagger algebra (resp. $k$-affinoid algebra) $A$. Then $d:A\to\Omega^1_X(X)=\Omega^1_A$ is the universal $k$-derivation of $A$ into finite $A$-modules. A dagger space $X$ is smooth if and only if its associated rigid space $X'$ is smooth, and $\Omega_{X'}^1$ is canonically identified with the coherent ${\cal O}_{X'}$-module obtained by completing on affinoid open pieces the coherent ${\cal O}_X$-module $\Omega_X^1$.\\
${\cal D}_X$ is a coherent ${\cal D}_X$-module (see \cite{doktor}, 5.5).\\If $X=\spm(A)$, then ${\cal D}_X(X)=D_A$ is a both-sided noetherian ring, and $\dim(A)=\kodi(D_A)$ (see \cite{mebdmod},p.42-55).\\

\addtocounter{satz}{1}{\bf \arabic{section}.\arabic{satz}}\newcounter{switchside1}\newcounter{switchside2}\setcounter{switchside1}{\value{section}}\setcounter{switchside2}{\value{satz}} Let $X$ be smooth of pure dimension. Then $\omega_X$ can be equipped with a canonical structure of ${\cal D}_X$-right-module, and there is an equivalence between the category of ${\cal D}_X$-right-modules and the category of ${\cal D}_X$-(left-)modules, 
where a ${\cal D}_X$-right-module ${\cal E}$ is sent to  ${\underline\ho}_{{\cal O}_X}(\omega_X,{\cal E})$, and a ${\cal D}_X$-(left-)module ${\cal M}$ is sent to $\omega_X\otimes_{{\cal O}_X}{\cal M}$. 
This can be seen as over the complex numbers, cf. \cite{mebdmod}.\\

\addtocounter{satz}{1}{\bf \arabic{section}.\arabic{satz}}  For a morphism $f:Z\to Y$ of smooth pure dimensional dagger (resp. rigid) spaces we define 
$${\cal D}_{Z\to Y}={\cal O}_Z\otimes_{f^{-1}{\cal O}_Y}f^{-1}{\cal D}_Y,$$ which is a $({\cal D}_Z, f^{-1}{\cal D}_Y)$-bimodule, and $${\cal D}_{Y\leftarrow Z}=\omega_Z\otimes_{f^{-1}{\cal O}_Y}     f^{-1}({\underline\ho}_{{\cal O}_Y}(\omega_Y,{\cal D}_Y))=\omega_{Z/Y}\otimes_{{\cal O}_Z}{\cal D}_{Z\to Y},$$ which is a $(f^{-1}{\cal D}_Y,{\cal D}_Z)$-bimodule, cf.\cite{mebdmod}.\\
In particular, if ${\cal M}$ is a ${\cal D}_Z$-module, ${\cal D}_{Y\leftarrow Z}\otimes_{{\cal D}_Z}{\cal M}$ becomes a $f^{-1}{\cal D}_Y$-module, and we get a left-derived functor $${\cal D}_{Y\leftarrow Z}\otimes_{{\cal D}_Z}^{\bf L}(.):\der^-({\cal D}_Z)\to \der^-(f^{-1}{\cal D}_Y).$$

\addtocounter{satz}{1}{\bf \arabic{section}.\arabic{satz}} \newcounter{spencer1}\newcounter{spencer2}\setcounter{spencer1}{\value{section}}\setcounter{spencer2}{\value{satz}} As in the complex case, we have a canonical ${\cal D}_X$-linear projective resolution \begin{gather} 0\to{\cal D}_X\otimes_{{\cal O}_X}\bigwedge^n{\cal T}_X\to\ldots\to{\cal D}_X\otimes_{{\cal O}_X}{\cal T}_X\to{\cal D}_X\tag{$*$}\end{gather} of ${\cal O}_X$, with ${\cal T}_X={\underline\ho}_{{\cal O}_X}(\Omega^1_X,{\cal O}_X)$. 
Therefore the definition $\dr({\cal M})=R{\underline\ho}_{{\cal D}_X}({\cal O}_X,{\cal M})$ makes sense for any ${\cal M}\in\der({\cal D}_X)$. 
For $r\in\mathbb{N}$ we define the $k$-vector spaces $H_{dR}^r(X)=H^r(X,\dr({\cal O}_X))$.\\
Application of $\omega_X\otimes_{{\cal O}_X}(.)$ to $(*)$ yields $\dr({\cal M})\cong\omega_X\otimes_{{\cal D}_X}^{\bf L}{\cal M}[-\dim(X)]$, and if $Y$ is another smooth space, $Z=X\times Y\to Y$ the canonical projection, there is a canonical isomorphism $\Omega^{\bullet}_{Z/Y}\cong{\cal D}_{Y\leftarrow Z}\otimes_{{\cal D}_Z}^{\bf L}{\cal O}_Z[-\dim(X)]$.\\

\addtocounter{satz}{1}{\bf \arabic{section}.\arabic{satz}} We recall a definition from \cite{huet},5.6. Let $f:Z\to Y$ be partially proper, let ${\cal F}$ be an abelian sheaf on $Z$. If $Y$ is quasi-compact, let $$\Gamma_c(Z/Y,{\cal F})=\{s\in\Gamma(Z,{\cal F}) |\mbox{ there is a quasi-compact}$$$$\mbox{admissible open } U\subset Z \mbox{ such that } s\in\ke(\Gamma(Z,{\cal F})\to\Gamma(Z-U,{\cal F}))\}.$$If $Y$ is arbitrary, let$$\Gamma_c(Z/Y,{\cal F})=\{s\in\Gamma(Z,{\cal F}) |\mbox{ for all quasi-compact}$$$$\mbox{admissible open } Y'\subset Y \mbox{ we have } s|_{f^{-1}(Y')}\in\Gamma_c(f^{-1}(Y')/Y',{\cal F}) \}.$$
Then $V\mapsto\Gamma_c(f^{-1}(V)/V,{\cal F})$ defines a sheaf $f_!{\cal F}$ on $Y$. We denote by $Rf_!:\der^+(Z)\to \der^+(Y)$ the functor induced by the left exact functor $f_!(-)$. In the following we always assume tacitly that $Rf_!:\der^+(Z)\to \der^+(Y)$ extends to $Rf_!:\der(Z)\to \der(Y)$ and $Rf_!:\der^-(Z)\to \der^-(Y)$ (this will be true in the cases relevant for us; namely, if $Z=X\times Y$ with $X$ finitely admissibly covered by open subspaces which are Zariski closed in some $({\bf D}^0)^m$, and $f$ is the projection --- as in \ref{spur} --- then there is an $n\in\mathbb{N}$ such that $R^if_!{\cal F}=0$ for all $i>n$, all ${\cal F}$. See \cite{doktor}, p.47).\\
If $Z$ and $Y$ are smooth and pure dimensional and if ${\cal M}\in\der({\cal D}_Z)$, we define $$f_+({\cal M})=Rf_!({\cal D}_{Y\leftarrow Z}\otimes_{{\cal D}_Z}^{\bf L}{\cal M})\in\der({\cal D}_Y).$$ If ${\cal M}\in\der^-({\cal D}_X)$ we get $f_+({\cal M})\in\der^-({\cal D}_Y)$.\\

\begin{pro} \label{kohkomp} Let $X\stackrel{f}{\to}Y\stackrel{g}{\to}Z$ be partially proper morphisms between smooth pure dimensional dagger  (resp. rigid) spaces. We assume that $g$ is a projection or a closed immersion. For ${\cal M}\in\der^-({\cal D}_X)$ there is a canonical isomorphism $g_+(f_+{\cal M})\cong(g\circ f)_+{\cal M}$ in $\der^-({\cal D}_Z)$. 
\end{pro}

{\sc Proof:} We have to show the projection formula$${\cal D}_{Z\leftarrow Y}\otimes_{{\cal D}_Y}^{\bf L}Rf_!({\cal D}_{Y\leftarrow X}\otimes_{{\cal D}_X}^{\bf L}{\cal M})\cong Rf_!(f^{-1}{\cal D}_{Z\leftarrow Y}\otimes_{f^{-1}{\cal D}_Y}^{\bf L}{\cal D}_{Y\leftarrow X}\otimes_{{\cal D}_X}^{\bf L}{\cal M}).$$If $g$ is a projection $Y=W\times Z\to Z$, we see that ${\cal D}_{Z\leftarrow Y}=\omega_{Y/Z}\otimes_{{\cal O}_Z}{\cal D}_Z$ is a coherent ${\cal D}_Y={\cal D}_W\otimes_k{\cal D}_Z$-right-module (because $\omega_{Y/Z}=\omega_W\otimes_k{\cal O}_Z$ is a coherent ${\cal D}_W\otimes_k{\cal O}_Z$-right-module). Therefore we obtain the above projection formula using the way-out-lemma (\cite{rd},I,7), since for finite free ${\cal D}_Y$-right-modules (instead of ${\cal D}_{Z\leftarrow Y}$) the projection formula is evident. 
If $g$ is a closed immersion, ${\cal D}_{Z\leftarrow Y}$ is a locally free ${\cal D}_Y$-right-module. This can be shown as in the complex analytic case using the fact that, locally for an admissible covering of $Z$, there are isomorphisms $Y\times{\bf D}^{\codim_Z(Y)}\cong Z$ such that $g$ corresponds to the zero section --- this is \cite{kidr},Theorem 1.18, in the rigid case, but holds true also in the dagger case as one observes by examining the proof in loc.cit.. 
By the commutation of $Rf_!$ with pseudo-filtered limits (cf. \cite{huet},5.3.7, or \cite{doktor}, 4.8) again one reduces the proof of the projection formula to the case of finite free ${\cal D}_Y$-right-modules.\\

\addtocounter{satz}{1}{\bf \arabic{section}.\arabic{satz}}  Let $Y\to X$ be a closed immersion into a smooth dagger (resp. rigid) space, defined by the coherent ideal ${\cal J}\subset{\cal
O}_X$. Then $${\underline
\Gamma}_{*Y}({\cal E})=\lim_{\stackrel{\to}{n}}{\underline\ho}_{{\cal O}_X}({\cal O}_X/{{\cal
J}^n},{\cal E})$$ and $${\cal E}(*Y)=\lim_{\stackrel{\to}{n}}{\underline\ho}_{{\cal O}_X}({{\cal
J}^n},{\cal E})$$ for a ${\cal D}_X$-module ${\cal E}$ are again ${\cal D}_X$-modules. We get right-derived functors $R{\underline
\Gamma}_{*Y}(-)$ and $R(-)(*Y)$ as functors $\der^+({\cal D}_X)\to\der^+({\cal D}_X)$, but also as functors $\der({\cal D}_X)\to\der({\cal D}_X)$ and $\der^-({\cal D}_X)\to\der^-({\cal D}_X)$. We have distinguished triangles $$R{\underline \Gamma}_{*Y}({\cal
K})\to{\cal K}\to R{\cal K}(*Y)\stackrel{+1}{\to}$$ for all ${\cal K}\in\der({\cal D}_X)$.\\

\addtocounter{satz}{1}{\bf \arabic{section}.\arabic{satz}} \newcounter{eigensch1}\newcounter{eigensch2}\setcounter{eigensch1}{\value{section}}\setcounter{eigensch2}{\value{satz}} We list some properties of the above functors. The proofs are similar to those in \cite{mebdmod}, the projection formulas needed can be justified as in \ref{kohkomp}.\\
(a) If $f:Z\to Y$ is a partially proper morphism between smooth pure dimensional dagger  (resp. rigid) spaces and if ${\cal M}\in\der^-({\cal D}_Z)$, there is a canonical isomorphism $\dr(f_+{\cal M})[\dim(Y)]\cong Rf_!\dr({\cal M})[\dim(Z)]$.\\(b) If in addition $T\to Y$ is a closed immersion and if $T_Z=T\times_YZ$, there is a canonical isomorphism $R{\underline
\Gamma}_{*T}(p_+{\cal M})\cong p_+(R{\underline\Gamma}_{*T_Z}{\cal M})$.\\(c) Let $X$ be smooth, let $Y_i\to X$ be closed immersions ($i=1,2$) and let ${\cal M}\in\der({\cal D}_X)$. 
There is a canonical isomorphism $$R{\underline\Gamma}_{*Y_1}(R{\underline\Gamma}_{*Y_2}({\cal M}))\cong R{\underline\Gamma}_{*(Y_1\cap Y_2)}({\cal M})$$ and a distinguished triangle $$R{\underline\Gamma}_{*(Y_1\cap Y_2)}({\cal M})\to R{\underline\Gamma}_{*Y_1}({\cal M})\oplus R{\underline\Gamma}_{*Y_2}({\cal M})\to R{\underline\Gamma}_{*(Y_1\cup Y_2)}({\cal M})\stackrel{+1}{\to}.$$(d) Let $X$ be smooth, let $Y\to X$ be a closed immersion and let ${\cal M}\in\der({\cal D}_X)$. 
There is a canonical isomorphism $R{\underline\Gamma}_{*Y}({\cal O}_X)\otimes_{{\cal O}_X}^{\bf L}{\cal M}\cong R{\underline\Gamma}_{*Y}({\cal M})$.\\(e) Let $X$ be smooth and affinoid, $\spm(B)=Y\to X$ a closed immersion of pure codimension $d$ such that all local rings $B_x$ (for $x\in\spm(B)=Y$) are locally complete intersections, and let ${\cal F}$ be a ${\cal D}_X$-module which is locally free as ${\cal O}_X$-module. 
Then $R^i{\underline\Gamma}_{*Y}({\cal F})=0$ for all $i\neq d$ (a well known algebraic fact!).\\

\begin{pro} \label{dmodgys} Let $Z\to Y\stackrel{s}{\to}X$ be a chain of closed immersions, where $Y$ and $X$ are smooth and pure dimensional. If $d=\codim(s)$, there is a canonical isomorphism $s_+R{\underline\Gamma}_{*Z}{\cal O}_Y\cong R{\underline\Gamma}_{*Z}{\cal O}_X[d]$.
\end{pro}

{\sc Proof:} (compare with \cite{mebend},3.3-1) First, \arabic{eigensch1}.\arabic{eigensch2}(b),(c) allows us to assume $Z=Y$. Denote by $I\subset{\cal O}_X$ the ideal of $Y$ in $X$. Observe ${\cal D}_{Y\to X}={\cal D}_X/I.{\cal D}_X$ and that this, as well as ${\cal D}_{X\leftarrow Y}$, is locally free over ${\cal D}_Y$, cf. the proof of \ref{kohkomp}. We claim that there is a canonical map of ${\cal D}_X$-right-modules$${\underline\ext}^d_{{\cal O}_X}({\cal O}_X/I,\omega_X)\otimes_{{\cal D}_Y}{\cal D}_{Y\to X}\to\lim_{\stackrel{\to}{k}}{\underline\ext}^d_{{\cal O}_X}({\cal O}_X/I^k,\omega_X).$$Indeed, choose an injective resolution $J^{\bullet}$ of the ${\cal D}_X$-right-module $\omega_X$ and define the morphism of complexes$${\underline\ho}_{{\cal O}_X}({\cal O}_X/I,J^{\bullet})\otimes_{{\cal D}_Y}{\cal D}_{Y\to X}\to\lim_{\stackrel{\to}{k}}{\underline\ho}_{{\cal O}_X}({\cal O}_X/I^k,J^{\bullet})$$as follows: If $g$ is a local section of ${\underline\ho}_{{\cal O}_X}({\cal O}_X/I,J^{m})$ and if $P$ is a local section of ${\cal D}_{Y\to X}$, represented by the local section $\tilde{P}$ of ${\cal D}_X$, then $g\otimes P$ is sent to the following local section of $\lim_{\stackrel{\to}{k}}{\underline\ho}_{{\cal O}_X}({\cal O}_X/I^k,J^{m})$: The ${\cal O}_X$-linear map ${\cal O}_X\to J^{m}$ which sends $1_{{\cal O}_X}$ to $g(1_{{\cal O}_X}).\tilde{P}$ actually induces an element of $\lim_{\stackrel{\to}{k}}{\underline\ho}_{{\cal O}_X}({\cal O}_X/I^k,J^{m})$: Indeed, if $\tilde{P}$ is of order $n$, then $g(1_{{\cal O}_X}).\tilde{P}$ is annihilated by $I^{n+1}$. We obtain the promised map. Now since we have as usual $\omega_Y\cong{\underline\ext}^d_{{\cal O}_X}({\cal O}_X/I,\omega_X)$, we get a map $$\omega_Y\otimes_{{\cal D}_Y}{\cal D}_{Y\to X}\to\lim_{\stackrel{\to}{k}}{\underline\ext}^d_{{\cal O}_X}({\cal O}_X/I^k,\omega_X)$$ of ${\cal D}_X$-right-modules. 
Applying ${\underline\ho}_{{\cal O}_X}(\omega_X,.)$ (cf. \arabic{switchside1}.\arabic{switchside2}) it becomes the map \begin{gather}{\cal D}_{X\leftarrow Y}\otimes_{{\cal D}_Y}{\cal O}_Y\to\lim_{\stackrel{\to}{k}}{\underline\ext}^d_{{\cal O}_X}({\cal O}_X/I^k,{\cal O}_X)=R^d{\underline\Gamma}_{*Y}{\cal O}_X\tag{$*$}\end{gather} of ${\cal D}_X$-left-modules. We claim that $(*)$ is an isomorphism. 
Indeed, if $x_1,\ldots, x_n$ are local coordinates on $X$ such that $Y$ is defined by $x_1,\ldots, x_d$ and if $\delta_1,\ldots,\delta_n$ is the basis of ${\underline\ho}_{{\cal O}_X}(\Omega_X,{\cal O}_X)$ dual to $dx_1,\ldots,dx_n$, one verifies that both sides in $(*)$ are identified with ${\cal D}_X/(x_1,\ldots,x_d,\delta_{d+1},\ldots,\delta_n)$. 
Since the right hand side in $(*)$ is already all of $R{\underline\Gamma}_{*Y}{\cal O}_X[d]$ (due to \arabic{eigensch1}.\arabic{eigensch2}(e)), and since on the left hand side we may write $\otimes^{\bf L}$ instead of $\otimes$, we are done.\\

\begin{pro}\label{spur}\label{kovhom}\label{funktkov} Let $X$ be smooth, proper and of pure dimension $n$, let $Y$ be smooth and pure dimensional and let $Z=X\times Y\stackrel{p}{\to}Y$ be the canonical projection. Then there is a canonical trace map $p_+{\cal O}_Z[n]\to{\cal O}_Y$. 
If furthermore $T\to Y$ and $S\to Z\times_YT=X\times T$ are closed immersions, there are canonical trace maps $p_+R{\underline\Gamma}_{*S}{\cal O}_Z[n]\to R{\underline\Gamma}_{*T}{\cal O}_Y$ and $$R\Gamma(Z,\dr(R{\underline\Gamma}_{*S}{\cal O}_Z))[2n]\to R\Gamma(Y,\dr(R{\underline\Gamma}_{*T}{\cal O}_Y)),$$ which are isomorphisms if the composition $S\to X\times T\to T$ is an isomorphism. 
\end{pro}

{\sc Proof:} We begin with the rigid case. In \cite{bey} it is described a finite admissible open covering $X=\cup_iU_i$ such that all $U_J=\cap_{i\in J}U_i$ (for $J\subset I$) have the following properties: 
there is a closed immersion $U_J\to({\bf D}^0)^{n_J}$ for some $n_J\in\mathbb{N}$, and for all affinoid $\hat{Y}$, all coherent ${\cal O}_{U_J\times\hat{Y}}$-modules ${\cal F}$ and all $j>n$ we have $R^j(p_{J,\hat{Y}})_!{\cal F}=0$, where $p_{J,\hat{Y}}:U_J\times\hat{Y}\to\hat{Y}$ denotes the projection. By means of Mayer-Vietoris sequences we get $R^jp_!{\cal F}=0$ for all coherent ${\cal O}_Z$-modules ${\cal F}$, all $j>n$, hence $R^mp_!\Omega^{\bullet}_{Z/Y}=0$ for $m>2n$ and $$R^{2n}p_!\Omega^{\bullet}_{Z/Y}=\frac{R^np_!\Omega^n_{Z/Y}}{\bi(R^np_!\Omega^{n-1}_{Z/Y}\to R^np_!\Omega^n_{Z/Y})}.$$ 
In view of $p_+{\cal O}_Z[n]\cong Rp_!\Omega^{\bullet}_{Z/Y}$ (cf. \arabic{spencer1}.\arabic{spencer2}), to give $p_+{\cal O}_Z[n]\to{\cal O}_Y$ it is therefore enough to give a map $t_p:R^np_!\omega_{Z/Y}\to{\cal O}_Y$ vanishing on $\bi(R^np_!\Omega^{n-1}_{Z/Y}\to R^np_!\Omega^n_{Z/Y})$. 
We take the following map (compare with \cite{vdpse}, \cite{bey}): For $$U_J\to D_J=({\bf D}^0)^{n_J}\subset{\bf D}^{n_J}=\spm(k<T_1,\ldots,T_{n_J}>)$$ as above and projection $p_{D_J,\hat{Y}}:D_J\times\hat{Y}\to\hat{Y}$, there is a canonical identification$$R^{n_J}(p_{D_J,\hat{Y}})_!\omega_{D_J\times\hat{Y}/\hat{Y}}(\hat{Y})=\{\omega=\sum_{\stackrel{\mu\in\mathbb{Z}^{n_J}}{\mu<0}}a_{\mu}T^{\mu}dT_1\wedge\ldots\wedge dT_{n_J}; a_{\mu}\in{\cal O}_{\hat{Y}}(\hat{Y})$$$$\mbox{ and }\omega\mbox{ converges on }\{t\in D_J;\mbox{ all }|t_i|<\epsilon\}\times\hat{Y}\mbox{ for some }0<\epsilon<1\}.$$For elements $\omega$ of this module set $t(\omega)=a_{(-1,\ldots,-1)}$. On the other hand, we have the Gysin map $g$:$$R^n(p_{J,\hat{Y}})_!\omega_{U_J\times\hat{Y}/\hat{Y}}(\hat{Y})\cong R^n(p_{J,\hat{Y}})_!\underline{Ext}^{n_J-n}({\cal O}_{U_J\times\hat{Y}},\omega _{D_J\times\hat{Y}/\hat{Y}})(\hat{Y})$$$$\to R^{n_J}(p_{D_J,\hat{Y}})_!\omega_{D_J\times\hat{Y}/\hat{Y}}(\hat{Y}).$$Locally it can be described as $\eta\mapsto \tilde{\eta}\wedge dX_{n+1}/X_{n+1}\wedge\ldots\wedge dX_{n_J}/X_{n_J}$ where $\tilde{\eta}$ is a lift of $\eta$ and $X_{n+1},\ldots, X_{n_J}$ are local equations for $U_J$ in $D_J$.\\We get $t\circ g:R^n(p_{J,\hat{Y}})_!\omega_{U_J\times\hat{Y}/\hat{Y}}(\hat{Y})\to{\cal O}_{\hat{Y}}(\hat{Y})$. This is seen to be independent of $n_J$ and of the chosen embedding $U_J\to({\bf D}^0)^{n_J}$, hence glues, for varying $J$, to give the desired map $t_p$. By construction, it vanishes on $\bi(R^np_!\Omega^{n-1}_{Z/Y}\to R^np_!\Omega^n_{Z/Y})$.\\In the dagger case, we argue by comparison with the associated morphism $p':Z'\to Y'$ of rigid spaces: Due to \cite{en1dag}, 3.5, we have again $R^jp_!{\cal F}=R^jp_*{\cal F}=0$ for all coherent ${\cal O}_Z$-modules ${\cal F}$, all $j>n$, and the composition of the canonical map $R^np_!\omega_{Z/Y}\to R^np'_!\omega_{Z'/Y'}$ with $t_{p'}:R^np'_!\omega_{Z'/Y'}\to{\cal O}_{Y'}$ has its image in ${\cal O}_Y\subset {\cal O}_{Y'}$, hence we obtain a map $t_p:R^np_!\omega_{Z/Y}\to{\cal O}_Y$ (this can be checked locally on $Y$; if $Y$ is affinoid, then this $t_p$ is the direct limit of the maps $t_{p_{\epsilon}}$ for the morphisms of rigid spaces $p_{\epsilon}:X'\times Y_{\epsilon}'\to Y_{\epsilon}'$, for appropriate extensions $Y'\subset Y_{\epsilon}'$).\\
If now in addition $T$ and $S$ are given, we can derive from $t_p$ the other promised maps using the isomorphisms from \arabic{eigensch1}.\arabic{eigensch2} (note that $$R\Gamma(Z,\dr(R{\underline\Gamma}_{*S}{\cal O}_Z))\cong R\Gamma(Y,Rp_!\dr(R{\underline\Gamma}_{*S}{\cal O}_Z))$$ because $S\to T$ is quasi-compact). 
Finally, suppose $S\to T$ is an isomorphism. Our additional statement in this situation is seen to be local on $X$. By the definition of $t_p:R^np_!\omega_{Z/Y}\to{\cal O}_Y$ we may substitute our $X$ by $X=({\bf D}^0)^n$ (dropping the assumption on properness). Passing to an admissible covering of $Y$, we may assume that there is a section $s:Y\to Z$ of $p:Z\to Y$ inducing the inverse of $S\to T$. It comes with an isomorphism$$R{\underline\Gamma}_{*T}{\cal O}_Y\cong R{\underline\Gamma}_{*T}p_+R{\underline\Gamma}_{*Y}{\cal O}_Z[n]\cong p_+R{\underline\Gamma}_{*S}{\cal O}_Z[n]$$by \arabic{eigensch1}.\arabic{eigensch2},\ref{dmodgys}. It is enough to show that its composition with the map in question $p_+R{\underline\Gamma}_{*S}{\cal O}_Z[n]\to R{\underline\Gamma}_{*T}{\cal O}_Y$ is an isomorphism. Of course this will follow once we know that the underlying map $${\cal O}_Y\cong p_+s_+{\cal O}_Y\cong p_+R{\underline\Gamma}_{*Y}{\cal O}_Z[n]\to p_+{\cal O}_Z[n]\to{\cal O}_Y$$is the identity. Since ${\cal O}_Y(Y)$ is jacobson, we may assume $Y=\spm(k)$ for this. The definition of $t_p$ as above does not depend on the choice of the closed embedding into some $({\bf D}^0)^n$: This tells us that the map $t_{id_Y}$ we get for $X=Y=\spm(k)$, computed by means of the embedding $s$, is the identity; but on the other hand, it is precisely the map we are interested in, by the compatibility of the Gysin map in the definition of $t_p$ with the Gysin map in \ref{dmodgys}.\\   

\begin{kor} \label{homdefbkt} Let $g_i:Z\to Y_i$ ($i=1,2$) be closed immersions into smooth affinoid dagger spaces $Y_i$ of pure dimension $n_i$. Then we have $$R\Gamma(Y_1,\dr(R{\underline\Gamma}_{*Z}{\cal O}_{Y_1}))[2n_1]\cong R\Gamma(Y_2,\dr(R{\underline\Gamma}_{*Z}{\cal O}_{Y_2}))[2n_2].$$
\end{kor}

{\sc Proof:} Because of \arabic{eigensch1}.\arabic{eigensch2},\ref{dmodgys} applied to closed immersions $Y_i\to{\bf D}^{n_i}$ we may assume $Y_i={\bf D}^{n_i}$. If $l:Z\to X=Y_1\times Y_2$ is the diagonal embedding it is enough to give isomorphisms $$\quad R\Gamma(Y_i,\dr(R{\underline\Gamma}_{*Z}{\cal O}_{Y_i}))\cong R\Gamma(X,\dr(R{\underline\Gamma}_{*Z}{\cal O}_X))[2n_{3-i}].$$ Let $i=1$. The open embedding $j:X\to W=Y_1\times{\bf P}^{n_2}_k$ induced by the open embedding into projective space $Y_2={\bf D}^{n_2}\to{\bf P}^{n_2}_k$ induces a closed immersion $(j\circ l):Z\to W$, and we have $R\Gamma(X,\dr(R{\underline\Gamma}_{*Z}{\cal O}_X))\cong R\Gamma(W,\dr(R{\underline\Gamma}_{*Z}{\cal O}_W))$. Now apply \ref{funktkov}.\\

\addtocounter{satz}{1}{\bf \arabic{section}.\arabic{satz}}  Let $Z$ be an affinoid $k$-dagger space, $q\in\mathbb{N}$. The definition $$h^{dR}_q(Z,k)=h^{dR}_q(Z)=\dim_k(H^{2n-q}(Y,\dr(R{\underline\Gamma}_{*Z}{\cal O}_Y))),$$ where $Z\to Y$ is a closed embedding into a smooth affinoid $k$-dagger space $Y$, is justified by \ref{homdefbkt}. For a finite field extension $k\subset k_1$ let $(?)_1=(?)\times_{\spm(k)}{\spm(k_1)}$. Then $h^{dR}_q(Z,k)=h^{dR}_q(Z_1,k_1)$. Indeed, clearly $\dim_k(H_{dR}^*(X/k))=\dim_{k_1}(H_{dR}^*(X_1/k_1))$ for any smooth $k$-dagger space $X$, hence $$\dim_k(H^*(X,\dr(R\underline{\Gamma}_Z{\cal O}_X)))=\dim_{k_1}(H^*(X_1,\dr(R\underline{\Gamma}_{Z_1}{\cal O}_{X_1})))$$ for closed subspaces $Z$ of smooth $k$-dagger spaces $X$. Now use \ref{algana} below.\\

\addtocounter{satz}{1}{\bf \arabic{section}.\arabic{satz}} \newcounter{endspur1}\newcounter{endspur2}\setcounter{endspur1}{\value{section}}\setcounter{endspur2}{\value{satz}}  Let $f:X\to Y$ be a finite \'{e}tale morphism of smooth dagger  (or rigid) spaces, $Y$ irreducible, $X=\cup X_i$ the decomposition into connected components. 
Assume all maps $f|_{X_i}:X_i\to Y$ to be surjective. Then there is an $l=\deg(f)\in\mathbb{N}$ and a trace map $$t:f_*\Omega^{\bullet}_X\to\Omega^{\bullet}_Y$$ such that the composition $\Omega_Y^{\bullet}\to f_*\Omega^{\bullet}_X\stackrel{t}{\to}\Omega^{\bullet}_Y$ is multiplication by $l$. 
In particular, $H^i_{dR}(Y)\to H^i_{dR}(X)$ is injective for all $i\in\mathbb{N}$.\\
Indeed, for admissible open connected $U=\spm(A)\subset Y$ with decomposition $f^{-1}(U)=\cup_j\spm(B_j)$ such that each $B_j$ is free over $A$, let $t_j:B_j\to A$ be the trace map, and for $q\in\mathbb{Z}$ let $f_*\Omega^q_X(X)\to\Omega^q_Y(Y)$ be the $A$-linear map $$f_*\Omega^q_X(X)=\Omega_{B_j}^q=\Omega_A^q\otimes_AB_j\to\Omega_A^q=\Omega^q_Y(Y)$$ which sends $\omega\otimes b$ to $t_j(b).\omega$. By the same computation as in \cite{hadr}, p.35, we see that for varying $q$ it commutes with the differentials. Clearly it glues for varying $U$, and the number $l=\sum_jl_j$, where $l_j$ denotes the rank of $B_j$ over $A$, is independent of $U$ and fulfills our requirement.\\  
(In fact, the \'{e}taleness of $f$ is not really needed: $f$ is flat in any case, by regularity of $X$ and $Y$. If $U$ is as above, let $L_j=\q(B_j)$, $K=\q(A)$; then $f$ induces finite separable field extensions $K\subset L_j$; let $l=\sum_j[L_j:K]$. The trace maps $\sigma_j:L_j\to K$ give rise to $\Omega_{L_j/k}^q=\Omega_{K/k}^q\otimes_KL_j\to\Omega_{K/k}^q, \omega\otimes f\mapsto \sigma_j(f).\omega$, restricting to $\sigma_j:\Omega_{B_j}^q\to\Omega_A^q$. Compare with the discussion in \cite{mw}, thm.8.3.. We do not need this.)\\

\addtocounter{satz}{1}{\bf \arabic{section}.\arabic{satz}}\newcounter{topsup1}\newcounter{topsup2}\setcounter{topsup1}{\value{section}}\setcounter{topsup2}{\value{satz}}  Let $X$ be a smooth dagger  (or rigid) space, $j:U\to X$ an open immersion with complement $Y=X-j(U)$. We do not put a structure of dagger (or rigid) space on $Y$. By $R{\underline\Gamma}_Y(.):\der^+({\cal D}_X)\to\der^+({\cal D}_X)$ we denote the right-derived functor of the left exact functor $${\cal F}\mapsto \ke({\cal F}\to j_*j^{-1}{\cal F})$$on abelian sheaves on $X$, and by $Rj_*:\der^+({\cal D}_U)\to\der^+({\cal D}_X)$ we denote the right-derived functor of $j_*$. 
Note that $Rj_*\dr({\cal L})\cong\dr(Rj_*{\cal L})$ for ${\cal L}\in\der^+({\cal D}_U)$. If $j':U'\to X$ is another open immersion with complement $Y'=X-j'(U')$ such that $U'\cup U$ is an admissible covering of an admissible open subset of $X$, there is a distinguished triangle $$R{\underline \Gamma}_{Y\cap Y'}({\cal K})\to
R{\underline \Gamma}_Y({\cal K})\oplus R{\underline
\Gamma}_{Y'}({\cal K})\to R{\underline \Gamma}_{Y\cup Y'}({\cal K})\stackrel{+1}{\to}.$$

\begin{pro}\label{algana}(a) Let $Y\to X$ be a closed immersion into a smooth dagger  (or rigid) space $X$. The canonical map \begin{gather}R\Gamma(X,\dr(R{\underline\Gamma}_{*Y}{\cal O}_X))\to R\Gamma(X,\dr(R{\underline\Gamma}_Y{\cal O}_X))\tag{$*$}\end{gather} is an isomorphism.\\(b) Let $Z\to Y$ be another closed immersion. 
There is a long exact sequence$$\ldots\to H^i(X,\dr(R{\underline\Gamma}_{*Z}{\cal O}_X))\to H^i(X,\dr(R{\underline\Gamma}_{*Y}{\cal O}_X))$$$$\to H^i(X-Z,\dr(R{\underline\Gamma}_{*Y-Z}{\cal O}_{X-Z}))\to H^{i+1}(X,\dr(R{\underline\Gamma}_{*Z}{\cal O}_X))\to\ldots.$$
\end{pro}

{\sc Proof:} {\bf (A):} The rigid case. First assume that $Y$ is locally defined by a single equation. Then the inclusion $j:U=(X-Y)\to X$ of the complement is a quasi-Stein morphism, hence (\cite{kiaub}) acyclic for coherent ${\cal O}_U$-modules (since we do not know if the analogue in the dagger case holds, we are forced to distinguish). It follows $\dr(Rj_*{\cal O}_U)=j_*\Omega^{\bullet}_U$. 
On the other hand $R{\cal O}_X(*Y)={\cal O}_X(*Y)$, and by \cite{kidr},Thm.2.3, the canonical map $DR({\cal O}_X(*Y))\to j_*\Omega^{\bullet}_U$ is an isomorphism. We get (a) for this type of $Y$. For general $Y$ assertion (a) is now deduced by an induction on the number of defining local equations, using the Mayer Vietoris sequences from \arabic{eigensch1}.\arabic{eigensch2} and \arabic{topsup1}.\arabic{topsup2}. Assertion (b) follows from (a) and the fact that for every sheaf ${\cal F}$ on $X$, we have a natural distinguished triangle $$R{\underline\Gamma}_Z{\cal F}\to R{\underline\Gamma}_Y{\cal F}\to Rj_*R{\underline\Gamma}_{Y-Z}j^{-1}{\cal F}\stackrel{+1}{\to}$$where $j:(X-Z)\to X$ is the open immersion: Take an injective resolution $I^{\bullet}$ of ${\cal F}$, then $j^{-1}I^{\bullet}$ is an injective resolution of $j^{-1}{\cal F}$, and $$0\to {\underline\Gamma}_Z I^{\bullet}\to {\underline\Gamma}_Y I^{\bullet}\to j_*{\underline\Gamma}_{Y-Z}j^{-1}I^{\bullet}\to 0$$ is exact.\\
{\bf (B):} The dagger case. Again (b) follows from (a). For (a), first assume that $Y$ is also smooth. As in the proof of \ref{kohkomp}, we find an affinoid admissible open covering $X=\cup_{i\in I}U_i$ such that for each $i\in I$ either $U_i\cap Y$ is empty or there exists an isomorphism $\phi_i:U_i\cong{\bf D}^m\times(U_i\cap Y)$ such that $U_i\cap Y\to U_i$ is the zero section. By a Cech argument one sees that it is enough to prove that for all finite and non-empty subsets $J$ of $I$, if we set $U_J=\cap_{i\in J}U_i$, the canonical map$$R\Gamma(U_J,\dr(R{\underline\Gamma}_{*Y}{\cal O}_X))\to R\Gamma(U_J,\dr(R{\underline\Gamma}_{Y}{\cal O}_X))$$is an isomorphism. If $U_J\cap Y$ is empty this is trivial, so we assume $U_J\cap Y$ is non-empty. Choose one $j\in J$. For $\epsilon\in|k^*|\cap]0,1]$ let ${\bf D}(\epsilon)$ be the closed disk of radius $\epsilon$ (with its dagger structure) and $$U_{j,J,\epsilon}=\phi_j^{-1}({\bf D}^m(\epsilon)\times(U_J\cap Y)).$$The set of the open subspaces $U_{j,J,\epsilon}$ is cofinal in the set of all open neighbourhoods of $U_J\cap Y$ in the affinoid space $U_{j,J,1}$. Since $U_J\cap U_{j,J,1}$ is such a neighbourhood, we find an $\epsilon_0$ such that $U_{j,J,\epsilon_0}\subset U_J$. Now observe that the canonical restriction maps$$R\Gamma(U_J,\dr(R{\underline\Gamma}_{*Y}{\cal O}_X))\to R\Gamma(U_{j,J,\epsilon_0},\dr(R{\underline\Gamma}_{*Y}{\cal O}_X))$$$$R\Gamma(U_J,\dr(R{\underline\Gamma}_{Y}{\cal O}_X))\to R\Gamma(U_{j,J,\epsilon_0},\dr(R{\underline\Gamma}_{Y}{\cal O}_X))$$are isomorphisms. Therefore we need to show that $$R\Gamma(U_{j,J,\epsilon_0},\dr(R{\underline\Gamma}_{*Y}{\cal O}_X))\to R\Gamma(U_{j,J,\epsilon_0},\dr(R{\underline\Gamma}_{Y}{\cal O}_X))$$is an isomorphism. In other words, we may assume from the beginning that $X={\bf D}^m\times Y$ and $Y\to X$ is the zero section. Let $D={\bf D}^m\subset{\bf P}^m_k=(\mbox{projective space})=P$, let $0$ be its origin, let $W=P\times Y$, $V=P-\{0\}$ and think of $Y=\{0\}\times Y$ as embedded into $W$. Since the natural restriction maps $$R\Gamma(W,\dr(R{\underline\Gamma}_{*Y}{\cal O}_W))\to R\Gamma(X,\dr(R{\underline\Gamma}_{*Y}{\cal O}_X))$$ $$R\Gamma(W,\dr(R{\underline\Gamma}_Y{\cal O}_W))\to R\Gamma(X,\dr(R{\underline\Gamma}_Y{\cal O}_X))$$are isomorphisms, it suffices to show that $$R\Gamma(W,\dr(R{\underline\Gamma}_{*Y}{\cal O}_W))\to R\Gamma(W,\dr(R{\underline\Gamma}_Y{\cal O}_W))$$ is an isomorphism. The dagger spaces $\{0\}, P$ and $V$ are partially proper, therefore {\bf (A)}(b) applies to give us the long exact Gysin sequence $$\ldots\to H_{dR}^{i-2m}(\{0\})\to
H_{dR}^i(P)\to H_{dR}^i(V)\to H_{dR}^{i-2m+1}(\{0\})\to\ldots.$$ By the K\"unneth formulas (in this case easily derived from \ref{semtub} below), we thus obtain the long exact sequence $$\ldots\to H_{dR}^{i-2m}(Y)\to
H_{dR}^i(W)\to H_{dR}^i(W-Y)\to H_{dR}^{i-2m+1}(Y)\to\ldots.$$Because of  $H_{dR}^{*-2m}(Y)\cong H^*(W,\dr(R{\underline\Gamma}_{*Y}{\cal O}_W))$, this implies what we want.\\Now for arbitrary $Y$, we may as in {\bf (A)} suppose that $Y$ is defined by a single equation and that $X$ is affinoid. 
Then we can reduce to the case where $Y$ is a divisor with normal crossings as in \cite{grodr}, considering a proper surjective morphism $g:X'\to X$ with $X'$ smooth, $U'=g^{-1}(U)\to U$ an isomorphism and $g^{-1}(Y)$ a divisor with normal crossings on $X'$ (such a $g$ exists by \ref{hironaka}). But in view of the Mayer Vietoris sequences from \arabic{eigensch1}.\arabic{eigensch2} and \arabic{topsup1}.\arabic{topsup2}, the normal crossings divisor case is equivalent with the case where $Y$ is smooth, which has been treated above.\\ 

\begin{kor} \label{blasexa} Let $f:X'\to X$ be a proper morphism of dagger  (or rigid) spaces, $Y\to X$ a closed immersion such that $f|_{X'-Y'}:(X'-Y')\to (X-Y)$ is an isomorphism, where $Y'=X'\times_XY$. 
Let $b:X\to Z$ be a closed immersion into a smooth space and $a:X'\to W'$ a locally closed immersion into a smooth proper space of pure dimension $n$. Then $(a,b\circ f):X'\to W'\times Z=Z'$ is a closed embedding, and there is a long exact sequence$$\ldots\to H^i(Z',\dr(R{\underline\Gamma}_{*Y'}{\cal O}_{Z'}))\to H^{i-2n}(Z,\dr(R{\underline\Gamma}_{*Y}{\cal O}_Z))\oplus H^i(Z',\dr(R{\underline\Gamma}_{*X'}{\cal O}_{Z'}))$$
$$\to H^{i-2n}(Z,\dr(R{\underline\Gamma}_{*X}{\cal O}_Z))\to H^{i+1}(Z',\dr(R{\underline\Gamma}_{*Y'}{\cal O}_{Z'}))\to\ldots.$$
\end{kor}

{\sc Proof:} By \ref{spur}, there is a morphism between the acyclic complexes which we get when we apply \ref{algana}(b) to $Y'\to X'\to Z'$ and to $Y\to X\to Z$; observe $$H^i(Z'-Y',\dr(R{\underline\Gamma}_{*X'}{\cal O}_{Z'}))\cong H^i(W'\times(Z-Y),\dr(R{\underline\Gamma}_{*X'}{\cal O}_{Z'}))$$for this, and that by the construction of the trace map in \ref{spur}, it is indeed a morphism of complexes, i.e. the resulting diagrams commute. Every third rung of this morphism of complexes is bijective (also by \ref{spur}), therefore we can perform a diagram chase according to the pattern \cite{hadr},p.44.\\ 

\begin{lem} \label{rigdrdf} Let $X_1$ and $X_2$ be smooth dagger spaces, let $\phi:X_1'\to X_2'$ be an isomorphism of the associated rigid spaces. Then $\phi$ gives rise to an isomorphism $\phi^{\dagger}:H^*_{dR}(X_2)\cong H^*_{dR}(X_1)$.
\end{lem}

{\sc Proof:} Set $X=X_1\times X_2$, $X'=X_1'\times X_2'$ and $\tilde{\phi}=(id,\phi):X_1'\to X'$, and let $\Delta=\bi(\tilde{\phi})$, a Zariski closed subspace of $X'$ isomorphic to $X_1'$. The canonical projections $X_1\leftarrow X\to X_2$ induce maps
$$H^*_{dR}(X_1)\stackrel{a_1}{\longrightarrow}H^*(X,\lim_{\stackrel{\to}{V}}j_{V*}\Omega^{\bullet}_V)\stackrel{a_2}{\longleftarrow}H^*_{dR}(X_2),$$where in the middle term $V$ runs through the open immersions $j_V:V\to X$ of dagger spaces with $\Delta\subset V'$, where $V'$ is the rigid space associated with $V$, regarded as an open subspace of $X'$. We claim that the $a_i$ are isomorphisms. The claim is local, so we may assume $X_1, X_2$ are affinoid and connected, there are elements $t_1,\ldots, t_m\in{\cal O}_{X_1}(X_1)=A_1$ such that $dt_1,\ldots, dt_m$ is a basis of $\Omega_{X_1}^1(X_1)$ over $A_1$, an open affinoid subspace $U\subset X$ such that $\Delta\subset U'$, where $U'\subset X'$ is the associated rigid space, an element $\delta\in\Gamma^*$, and an isomorphism of rigid spaces$$\rho:U'\to\spm(k<\delta^{-1}.Z_1,\ldots,\delta^{-1}.Z_m>)\times\Delta$$
where $\delta^{-1}.Z_i$ is sent to $\delta^{-1}.(t_i\otimes1-1\otimes(\phi^*)^{-1}(t_i))\in {\cal O}_{U'}(U')$, compare \cite{kidr},1.18. 
For $0<\epsilon\le\delta$ let $$U_{\epsilon}'=\rho^{-1}(\spm(k<\epsilon^{-1}.Z_1,\ldots,\epsilon^{-1}.Z_m>)\times\Delta),$$ an open subspace of $X'$, and let $U_{\epsilon}$ be the corresponding open subspace of $X$. Since $U_{\epsilon}'$ is a Weierstrass domain in $X'$, the same is true for $U_{\epsilon}$ in $X$ (if necessary, modify the defining functions slightly to get overconvergent ones); in particular, $U_{\epsilon}$ is affinoid, so $j_{{U_{\epsilon}}*}\Omega^{\bullet}_{U_{\epsilon}}=Rj_{{U_{\epsilon}}*}\Omega^{\bullet}_{U_{\epsilon}}$. Since $X$ is quasi-compact, $\lim_{\stackrel{\to}{V}}$ commutes with formation of cohomology, and since the ${U_{\epsilon}}$ are cofinal in $\{V\}$, it is now enough to show that for arbitrary $0<\epsilon\le\delta$ the maps $$b_{i,\epsilon}:H_{dR}^*(X_i)\to H_{dR}^*(U_{\epsilon})$$ are isomorphisms ($i=1,2$). 
By \cite{boartin}, we can find an isomorphism $\sigma:X_1\to X_2$ such that the induced map $\sigma':X_1'\to X_2'$ is close to $\phi$; in particular so close that for $\tilde{\sigma}=(id,\sigma):X_1\to X$ we have $\bi(\tilde{\sigma})\subset U_{\epsilon/2}$. Similarly, we can approximate the map ${\cal O}_{X_1'}(X_1')\cong{\cal O}_{\Delta}(\Delta)\stackrel{\rho^*}{\to}{\cal O}_{U'}(U')\to{\cal O}_{U_{\epsilon}'}(U_{\epsilon}')$ by a map ${\cal O}_{X_1}(X_1)\to{\cal O}_{U_{\epsilon}}(U_{\epsilon})$. 
Its extension to the map$${\cal O}_{X_1}(X_1)\otimes_k^{\dagger}k<\epsilon^{-1}.Z_1,\ldots,\epsilon^{-1}.Z_m>^{\dagger}\to{\cal O}_{U_{\epsilon}}(U_{\epsilon})$$which sends $\epsilon^{-1}.Z_i$ to  $\epsilon^{-1}.(t_i\otimes1-1\otimes(\sigma^*)^{-1}(t_i))$ is an isomorphism since its completion is close to the isomorphism obtained from $\rho$. 
So we have an isomorphism $$U_{\epsilon}\cong\spm(k<\epsilon^{-1}.Z_1,\ldots,\epsilon^{-1}.Z_m>^{\dagger})\times X_1$$ where the closed immersion $\tilde{\sigma}:X_1\to U_{\epsilon}$ corresponds to the zero section. 
Hence the maps $H_{dR}^*(U_{\epsilon})\to H_{dR}^*(X_1)$ induced by $\tilde{\sigma}$ are isomorphisms, by \ref{semtub}. Since $\tilde{\sigma}$ is a section for the canonical map $U_{\epsilon}\to X_1$ which gives rise to $b_{1,\epsilon}$, we derive the bijectivity of $b_{1,\epsilon}$. That $b_{2,\epsilon}$ is bijective is seen symmertrically. Now we define $\phi^{\dagger}=a_1^{-1}\circ a_2$.\\One can show that this construction is compatible with compositions: If $X_3$ is a third dagger space with associated rigid space $X_3'$, and if $\gamma:X_2'\to X_3'$ is an isomorphism, than $\phi^{\dagger}\circ\gamma^{\dagger}=(\gamma\circ\phi)^{\dagger}$, see \cite{en3dag}. We do not need this here.\\

\section{De Rham cohomology of tubes of a semi-stable reduction}
From now on let $R$ be a complete discrete valuation ring of mixed characteristic $(0,p)$, let $\pi\in R$ be a uniformizer, $k$ its fraction field, and $\bar{k}=R/(\pi)$ its residue field.

\begin{lem}\label{semtub} (a) Let $r, n\in\mathbb{Z}$, $0\le r\le n$, let $\mu$ and $\gamma_i$ for $1\le i\le r$, and $\delta_i$ for $1\le i\le n$ be elements of $\Gamma^*$ such that $\gamma_i\le\delta_i$ for all $1\le i\le r$, and $\prod_{i=1}^r\delta_i\ge \mu$. Define the open dagger subspace $V$ of the dagger affine space ${\bf A}_k^n$ by$$V=\{(x_1,\ldots,x_n)\in{\bf A}_k^n|$$$$\prod_{i=1}^r|x_i|\ge\mu, |x_i|\le\delta_i \mbox{ for all }1\le i\le n,\mbox{ and }|x_i|\ge\gamma_i \mbox{ for all }1\le i\le r\}.$$Let $X_1,\ldots, X_r$ be the first $r$ standard coordinates on ${\bf A}_k^n$. Then $H_{dR}^q(V)$ is the $k$-vector space generated by the classes of the $q$-forms $$dX_{i_1}/X_{i_1}\wedge\ldots\wedge dX_{i_q}/X_{i_q}$$ with $1\le i_1<\ldots<i_q\le r$. In particular, if $r=0$ we have $H_{dR}^q(V)=0$ for all $q>0$. If $X$ is another smooth dagger space, the canonical maps$$\oplus_{q_1+q_2=q}H_{dR}^{q_1}(X)\otimes_kH_{dR}^{q_2}(V)\stackrel{\beta}{\to}H_{dR}^q(X\times V)$$ are bijective.\\
(b) Suppose even $\gamma_i<\delta_i$ for all $1\le i\le r$, and $\prod_{i=1}^r\delta_i<\mu$. Define the open dagger (resp. rigid) subspace $V$ of the dagger (resp. rigid) affine space ${\bf A}_k^n$ by$$V=\{(x_1,\ldots,x_n)\in{\bf A}_k^n|$$$$\prod_{i=1}^r|x_i|>\mu, |x_i|<\delta_i \mbox{ for all }1\le i\le n,\mbox{ and }|x_i|>\gamma_i \mbox{ for all }1\le i\le r\}.$$ Then the same assertions as in (a) hold (of course, if $V$ is the dagger (resp. rigid) space, then $X$ should be a dagger (resp. rigid) space, too). 
\end{lem}

{\sc Proof:} (a) Note that $V$ is affinoid. We may assume that $X$ is also affinoid and connected, $X=\spm(B)$. After a finite extension of $k$ we may assume there are $\underline{\delta}_i, \underline{\gamma}_i$ and $\underline{\mu}$ in $k$ such that $|\underline{\delta}_i|=\delta_i$, $|\underline{\gamma}_i|=\gamma_i$ and $|\underline{\mu}|=\mu$. We regard ${\cal O}_{X\times V}(X\times V)$ as a subring of $$B<\underline{\delta}_1^{-1}.X_1,(\underline{\delta}_1^{-1}.X_1)^{-1},\ldots,\underline{\delta}_n^{-1}.X_n,(\underline{\delta}_n^{-1}.X_n)^{-1}>$$(to compute in this ring was suggested by the referee).\\(i) We begin with the following observation. Let $$a=\sum_{v\in\mathbb{Z}^n}\theta_v\prod_{j=1}^nX_j^{v_j}\in{\cal O}_{X\times V}(X\times V),$$$\theta_v\in B$, $\theta_v=0$ whenever there is a $r< j\le n$ with $v_j<0$. 
Fix $1\le l\le n$. We claim$$b=\sum_{\stackrel{v\in\mathbb{Z}^n}{v_l\ne0}}v_l^{-1}\theta_v\prod_{j=1}^nX_j^{v_j}\in{\cal O}_{X\times V}(X\times V),$$i.e. also this sum converges in ${\cal O}_{X\times V}(X\times V)$. Indeed, consider the surjection of dagger algebras
$$D=B\otimes_k^{\dagger}k<X_1,\ldots,X_n,Y_1,\ldots,Y_r,Z>^{\dagger}\stackrel{\tau}{\longrightarrow}{\cal O}_{X\times V}(X\times V),$$$$X_i\mapsto\underline{\delta}_i^{-1}.X_i,\quad\quad Y_i\mapsto(\underline{\gamma}_i^{-1}.X_i)^{-1},\quad\quad Z\mapsto(\underline{\mu}^{-1}.\prod_{i=1}^rX_i)^{-1}.$$By definition of $\otimes_k^{\dagger}$, we have
$$D=\lim_{\stackrel{\to}{X_{\epsilon},\delta}}{\cal O}_{X_{\epsilon}}(X_{\epsilon})\hat{\otimes}_kk<\delta^{-1}.X_1,\ldots,\delta^{-1}.X_n,\delta^{-1}.Y_1,\ldots,\delta^{-1}.Y_r,\delta^{-1}.Z>$$where the $X_{\epsilon}$ run through the strict neighbourhoods of $X'$, the rigid space associated with $X$, in an appropriate affinoid rigid space which contains $X'$ as a relatively compact open subset, and where $\delta$ runs through all $\delta>1$. 
Hence, a given  $$c=\sum_{w\in\mathbb{Z}_{\ge0}^{n+m+1}}\beta_w(\prod_{j=1}^nX_j^{w_j})(\prod_{j=1}^rY_j^{w_{j+n}})Z^{w_{n+r+1}}\in D$$$(\beta_w\in B)$ is an element of $${\cal O}_{X_{\epsilon}}(X_{\epsilon})\hat{\otimes}_kk<\delta^{-1}.X_1,\ldots,\delta^{-1}.X_n,\delta^{-1}.Y_1,\ldots,\delta^{-1}.Y_r,\delta^{-1}.Z>$$ for some $X_{\epsilon},\delta$, and one easily sees that $$d=\sum_{\stackrel{w\in\mathbb{Z}_{\ge0}^{n+m+1}}{w_l-w_{l+n}-w_{n+r+1}\ne 0}}(w_l-w_{l+n}-w_{n+r+1})^{-1}\beta_w(\prod_{j=1}^nX_j^{w_j})(\prod_{j=1}^rY_j^{w_{j+n}})Z^{w_{n+r+1}}$$
is then an element of $${\cal O}_{X_{\epsilon}}(X_{\epsilon})\hat{\otimes}_kk<\delta_1^{-1}.X_1,\ldots,\delta_1^{-1}.X_n,\delta_1^{-1}.Y_1,\ldots,\delta_1^{-1}.Y_r,\delta_1^{-1}.Z>$$ for any $1<\delta_1<\delta$; in particular $d\in D$, too. Clearly, if $\tau(c)=a$, then $\tau(d)=b$, and the claim follows.\\ 
(ii) For $1\le i_1<\ldots<i_t\le n$ we write $$dX_{i}/X_i=dX_{i_1}/X_{i_1}\wedge\ldots\wedge dX_{i_t}/X_{i_t}$$
$$(dX_{i}/X_i)^{\hat{i_t}}=dX_{i_1}/X_{i_1}\wedge\ldots\wedge dX_{i_{t-1}}/X_{i_{t-1}}.$$
Every $\omega\in\Omega^q_{X\times V}(X\times V)$ can uniquely be written as a convergent series \begin{gather} \omega=\sum_{0\le t\le q}\sum_{1\le i_1<\ldots<i_t\le n}\sum_{v\in\mathbb{Z}^n}\sigma_{t,i,v}(\prod_{j=1}^nX_j^{v_j})dX_{i}/X_i\tag{$*$}\end{gather} with $\sigma_{t,i,v}\in\Omega^{q-t}_X(X)$ and $\sigma_{t,i,v}=0$ whenever there is a $j$ with $r<j\le n$ and $v_j\le 0$. For $0\le l\le n$ let $\Omega_l^q\subset \Omega^q_{X\times V}(X\times V)$ be defined by the additional condition: 
$\sigma_{t,i,v}=0$ whenever $l<i_t$ or $v_j\ne0$ for some $j$ with $l<j\le n$; this condition means, that no $X_j$, no $dX_j/X_j$ for $l<j$ occurs in $\omega$. We claim that, if $l>0$, every $\omega\in\Omega^q_l$ with $d\omega=0$ can be written modulo exact forms as $\omega=\omega_0+\omega_1dX_l/X_l$ with $\omega_0\in\Omega^q_{l-1}$ and $\omega_1\in\Omega^{q-1}_{l-1}$ and $d\omega_0=0$ and $d\omega_1=0$, and even $\omega_1=0$ if $l>r$. 
Indeed, if $\omega$ is represented as in $(*)$, then by (i) the series$$\eta=\sum_{0\le t\le q}\sum_{1\le i_1<\ldots<i_t=l}\sum_{\stackrel{v\in\mathbb{Z}^n}{v_l\ne0}}v_l^{-1}\sigma_{t,i,v}(\prod_{j=1}^nX_j^{v_j})(dX_{i}/X_i)^{\hat{i_t}}$$converges in $\Omega^{q-1}_l$, and subtracting $d\eta$ from $\omega$ we see that we may suppose $\sigma_{t,i,v}=0$ whenever $i_t=l$ and $v_l\ne0$. But $d\omega=0$ implies $\sigma_{t,i,v}=0$ whenever $v_l\ne0$ and $i_t<l$, so in any case we have $\sigma_{t,i,v}=0$ if only $v_l\ne0$. From this the claim follows.\\
(iii) Let $\omega\in\Omega^q_{X\times V}(X\times V)=\Omega_n^q$ be such that $d\omega=0$. Repeated application of (ii) shows that modulo exact forms $\omega$ can be written as $$\omega=\sum_{0\le t\le q}\sum_{1\le i_1<\ldots<i_t\le r}\sigma_{t,i}dX_{i}/X_i$$with uniquely determined $\sigma_{t,i}\in\Omega^{q-t}_X(X)$ such that $d\sigma_{t,i}=0$. This provides us with an inverse map for $\beta$ and proves the assertion on $H_{dR}^q(V)$.\\(b) Again we may assume $X$ is affinoid. 
In the dagger case, we then exhaust $V$ by affinoid dagger spaces of the type considered in (a) and conclude by using our result in (a) (passing to the limit). In the rigid case, $X\times V$ is quasi-Stein, hence acyclic for coherent modules (\cite{kiaub}). Therefore, in this case, we can argue literally as in (a) --- since here $V$ is defined by strict inequalities, we do not need overconvergence.\\

\addtocounter{satz}{1}{\bf \arabic{section}.\arabic{satz}} \newcounter{defsemist1}\newcounter{defsemist2}\setcounter{defsemist1}{\value{section}}\setcounter{defsemist2}{\value{satz}} We call a closed immersion ${\cal Z}\to{\cal X}$ of noetherian $\pi$-adic formal $R$-schemes a strictly semi-stable formal pair $({\cal X},{\cal Z})$ over $R$, if there is an $n\in\mathbb{N}$, a Zariski open covering ${\cal X}=\cup_i{\cal U}_i$ and for all $i$ a pair $s(i), r(i)\in\mathbb{N}$ with $n\ge s(i)\ge r(i)\ge 1$ and an \'{e}tale morphism
$$q_i:{\cal U}_i\to\spf(R<X_1,\ldots,X_n>/(X_1\ldots X_{r(i)}-\pi))$$such that ${\cal Z}|_{{\cal U}_i}=\cup_{j=r(i)+1,\ldots,s(i)}V(q_i^*X_j)$. We call ${\cal X}$ a strictly semi-stable formal $R$-scheme if $({\cal X},\emptyset)$ is a strictly semi-stable formal pair over $R$.\\

\addtocounter{satz}{1}{\bf \arabic{section}.\arabic{satz}} \newcounter{schwfmsch1}\newcounter{schwfmsch2}\setcounter{schwfmsch1}{\value{section}}\setcounter{schwfmsch2}{\value{satz}} For a $\pi$-adic topologically finite type (tf) formal $R$-scheme ${\cal X}$ with generic fibre (\cite{bolu}) the rigid space ${\cal X}_k$, there is a specialization map $s:{\cal X}_k\to{\cal X}$, and if $Y\to{\cal X}_{\bar{k}}$ is an immersion into the closed fibre, then $]Y[_{\cal X}=s^{-1}(X)$ is an admissible open subspace of ${\cal X}_k$, the tube of $Y$. 

\begin{satz}\label{monovor} Let ${\cal X}$ be a strictly semi-stable formal $R$-scheme, let ${\cal X}_{\bar{k}}=\cup_{i\in I}Y_i$ be the decomposition of the closed fibre into irreducible components. For $K\subset I$ set $Y_K=\cap_{i\in K}Y_i$. Let $X^{\dagger}$ be a $k$-dagger space such that its associated rigid space is identified with ${\cal X}_k$. For a subscheme $Y\subset{\cal X}_{\bar{k}}$ let $]Y[_{\cal X}^{\dagger}$ be the open dagger subspace of $X^{\dagger}$ corresponding to the open rigid subspace $]Y[_{\cal X}$ of ${\cal X}_k$. Then for any $\emptyset\ne J\subset I$ the canonical map $$H_{dR}^*(]Y_J[_{\cal X}^{\dagger})\to H_{dR}^*(]Y_J-(Y_J\cap(\cup_{i\in I-J}Y_i))[^{\dagger}_{\cal X})$$ is bijective.
\end{satz}
 
{\sc Proof:} (i) Suppose $I\ne J$. For $L\subset(I-J)$ set $$G_L=]Y_J-(Y_J\cap(\cup_{i\in L}Y_i))[^{\dagger}_{\cal X}.$$ For $m\in\mathbb{N}$ let $P_m(I-J)$ be the set of subsets of  $I-J$ with $m$ elements, and $G^m=\cup_{L\in P_m(I-J)}G_L$. Let $G^0=]Y_J[^{\dagger}_{\cal X}$. Then $G^{m+1}\subset G^m$ for all $m\ge0$, and $$G^{|I-J|}=]Y_J-(Y_J\cap(\cup_{i\in I-J}Y_i))[^{\dagger}_{\cal X}.$$ It is enough to show that \begin{gather}H_{dR}^*(G^m)\to H_{dR}^*(G^{m+1})\tag*{$(*)_m$}\end{gather} is bijective for all $m\le|I-J|-1$. Since $G^m=\cup_{L\in P_m(I-J)}G_L$ is an admissible covering, it is enough to show that \begin{gather} H_{dR}^*(\cap_{L\in Q}G_L)\to H_{dR}^*(\cap_{L\in Q}G_L\cap G^{m+1})\tag*{$(**)_m$}\end{gather} is bijective for all $Q\subset P_m(I-J)$. 
But $\cap_{L\in Q}G_L=G_M$ for $M=\cup_{L\in Q}L$, and $$\cap_{L\in Q}G_L\cap G^{m+1}=\cup_{i\in(I-(M\cup J))}(G_M-(G_M\cap]Y_i[^{\dagger}_{\cal X})).$$ If now ${\cal X}$ is replaced by its open formal subscheme which on the underlying topological space is the complement of $\cup_{i\in M}Y_i$ in ${\cal X}_{\bar {k}}$, then this means the replacement of $]Y_J[^{\dagger}_{\cal X}$ by $G_M$ and of $I$ by $I-M$ (but $J$ stays the same). In this way $(**)_m$ takes the form $(*)_0$, therefore it suffices to prove $(*)_0$. Note $G^1=]Y_J[^{\dagger}_{\cal X}-]Y_I[^{\dagger}_{\cal X}$, i.e. we must prove that $$H_{dR}^*(]Y_J[_{\cal X}^{\dagger})\to H_{dR}^*(]Y_J-Y_I[^{\dagger}_{\cal X})$$ is bijective.\\
(ii) Let ${\cal X}=\cup_{s\in S}{\cal U}_s$ be an open covering of ${\cal X}$ as in the definition of strict semi-stability. For a finite non-empty subset $T$ of $S$ let ${\cal U}_T=\cap_{s\in T}{\cal U}_s$. It is enough to show that for each such $T$ the map $$H_{dR}^*(]Y_J\cap{\cal U}_T[_{{\cal U}_T}^{\dagger})\to H_{dR}^*(](Y_J-Y_I)\cap{\cal U}_T[_{{\cal U}_T}^{\dagger})$$is bijective. This is trivial if $Y_I\cap{\cal U}_T$ is empty. If $Y_I\cap{\cal U}_T$ is non-empty, the irreducible components of the reduction $({\cal U}_T)_{\bar{k}}$ of ${\cal U}_T$ correspond bijectively to those of ${\cal X}_{\bar{k}}$, so we can replace ${\cal X}$ by ${\cal U}_T$. In other words, it is enough to prove the statement in (i) in the following case: ${\cal X}=\spf(A)$ is affine, and there is an \'{e}tale morphism ${\cal X}=\spf(A)\stackrel{\phi}{\to}\spf(R<X_1,\ldots,X_n>/(X_1\ldots X_r-\pi))$. Let $f_i=\phi^*(X_i)\in A$, inducing $\bar{f}_i\in A/(\pi)$. Passing to an open covering of ${\cal X}$ we may suppose that each $V(\bar{f_i})$ is irreducible (and non-empty), so we derive an identification $\{1,\ldots,r\}=I$. For $\lambda\in\Gamma^*$ with $\lambda<1$ set $$F_{\lambda}=\{x\in G^0|\quad|f_j(x)|\le\lambda\mbox{ for all }j\in J\}$$ and $E_{\lambda}=F_{\lambda}\cap G^1$. Then $G^0=\cup_{\lambda<1}F_{\lambda}$ is an admissible covering, and to prove $(*)_0$ it suffices to show the bijectivity of $H_{dR}^*(F_{\lambda})\to H_{dR}^*(E_{\lambda})$ for all such $\lambda$.\\
(iii) For $\beta\in\Gamma^*$ with $\beta<1$ and $i\in I-J$ set $$F_{\lambda,\beta}^i=\{x\in F_{\lambda}|\quad|f_i(x)|\ge\beta\}\quad\mbox{  and  }\quad F_{\lambda,\beta}=\cup_{i\in I-J}F_{\lambda,\beta}^i.$$ We have $E_{\lambda}\subset F_{\lambda,\beta}\subset F_{\lambda}$, and it suffices now to prove that the following maps are bijective:\begin{gather}\quad\lim_{\stackrel{\to}{\beta\to 1}}H^*_{dR}(F_{\lambda,\beta})\to H^*_{dR}(E_{\lambda})\tag{1}\\\quad H^*_{dR}(F_{\lambda})\to H^*_{dR}(F_{\lambda,\beta})\quad\mbox{  for }\beta<1.\tag{2}\end{gather} Note $G_{\{i\}}=]Y_J-(Y_J\cap Y_i)[^{\dagger}_{\cal X}=\{x\in]Y_J[^{\dagger}_{\cal X}; |f_i(x)|=1\}$ for $i\in I-J$. We compare the admissible covering $E_{\lambda}=\cup_{i\in I-J}(G_{\{i\}}\cap F_{\lambda})$ with the admissible covering $F_{\lambda,\beta}=\cup_{i\in I-J}F_{\lambda,\beta}^i$: since the direct limit is exact, to prove the bijectivity of (1), it is enough to prove the bijectivity of \begin{gather}\quad\lim_{\stackrel{\to}{\beta\to 1}}H^*_{dR}(\cap_{i\in K}F^i_{\lambda,\beta})\to H^*_{dR}(\cap_{i\in K}G_{\{i\}}\cap F_{\lambda})\tag{3}\end{gather}for all $\emptyset\ne K\subset(I-J)$. But the affinoid dagger space $\cap_{i\in K}G_{\{i\}}\cap F_{\lambda}$ is the inverse limit of the affinoid dagger spaces $\cap_{i\in K}F^i_{\lambda,\beta}$, in particular $$\Gamma(\cap_{i\in K}G_{\{i\}}\cap F_{\lambda},\Omega^{\bullet})=\lim_{\stackrel{\to}{\beta\to 1}}\Gamma(\cap_{i\in K}F^i_{\lambda,\beta},\Omega^{\bullet}),$$ so the bijectivity of (3) follows from the exactness of direct limits.\\(iv) It remains to prove the bijectivity of the maps (2). Set $$S_{\lambda,\beta}=\{x\in F_{\lambda}|\quad|f_i(x)|\le\beta\mbox{ for all }i\in I-J\}.$$ Then $F_{\lambda}=S_{\lambda,\beta}\cup F_{\lambda,\beta}$ is an admissible covering, and the bijectivity of (2) is equivalent to that of \begin{gather}\quad H^*_{dR}(S_{\lambda,\beta})\to H^*_{dR}(S_{\lambda,\beta}\cap F_{\lambda,\beta}).\tag{4}\end{gather}
(v) We claim that there is a $\pi$-adic affine formally smooth tf formal $\spf(R)$-scheme $\spf(D)$ and an isomorphism of rigid spaces\newpage $$]Y_I[_{{\cal X}}=]V((\bar{f}_i)_{i=1,\ldots,r})[_{{\cal X}}\cong$$$$\spm(D\otimes_Rk)\times\{x\in\spm(k<T_1,\ldots,T_r>/(T_1\ldots T_r-\pi))|\quad|T_i(x)|<1\mbox{ for all }i\}$$such that $f_i\in A$ corresponds to $T_i$. 
This is constructed as follows: Let $\hat{A}$ be the $(f_1,\ldots,f_r)$-adic completion of $A$. Since $\hat{A}/(\pi,f_1,\ldots,f_r)=A/(\pi,f_1,\ldots,f_r)=\bar{D}$ is a smooth $\bar{k}$-algebra ($\phi$ is \'{e}tale), there is a section $\bar{D}\stackrel{\bar{s}}{\to}\hat{A}/(\pi)$ for the canonical surjection $\hat{A}/(\pi)\to\bar{D}$, a lift of $\bar{D}$ to a smooth $R$-algebra $\tilde{D}$ (see \cite{elk}) and a lift of $\bar{s}$ to a morphism $D\stackrel{s}{\to}\hat{A}$, where $D$ is the $\pi$-adic completion of $\tilde{D}$. The extension $D[[T_1,\ldots,T_r]]\to\hat{A},\quad T_i\mapsto f_i$ of $s$ induces an isomorphism $$D[[T_1,\ldots,T_r]]/(T_1\ldots T_r-\pi)\cong\hat{A}.$$ This gives the desired isomorphism of rigid spaces (compare \cite{berco}, 0.2.7, for the computation of tubes).\\ 
(vi) To prove the bijectivity of (4), we may now, in view of \ref{rigdrdf}, assume that there is a smooth $k$-dagger algebra $B$ and an isomorphism of dagger spaces $$]Y_I[^{\dagger}_{{\cal X}}=]V((\bar{f}_i)_{i=1,\ldots,r})[^{\dagger}_{{\cal X}}\cong$$$$\spm(B)\times\{x\in\spm(k<T_1,\ldots,T_r>^{\dagger}/(T_1\ldots T_r-\pi))|\quad|T_i(x)|<1\mbox{ for all }i\}$$such that $f_i\in A$ corresponds to $T_i$. Let $$N=\{x\in S_{\lambda,\beta}|\quad|f_i(x)|=\beta\mbox{ for all }i\in I-J\}.$$ It suffices to show the bijectivity of the two maps \begin{gather}\quad H^*_{dR}(S_{\lambda,\beta})\to H^*_{dR}(N)\tag{5}\\\quad H^*_{dR}(S_{\lambda,\beta}\cap F_{\lambda,\beta})\to H^*_{dR}(N).\tag{6}\end{gather}The bijectivity of $(5)$ follows immediately from \ref{semtub}. Finally, consider the admissible covering $S_{\lambda,\beta}\cap F_{\lambda,\beta}=\cup_{i\in I-J}(S_{\lambda,\beta}\cap F_{\lambda,\beta}^i)$. To prove the bijectivity of $(6)$, it is enough to prove that of \begin{gather} H^*_{dR}(\cap_{i\in K}S_{\lambda,\beta}\cap F^i_{\lambda,\beta})\to H^*_{dR}(N)\tag{7}\end{gather}for all $\emptyset\ne K\subset(I-J)$, which again can be done using \ref{semtub}.  

\begin{satz}\label{maxschni} In \ref{monovor}, suppose in addition that ${\cal X}$ is quasi-compact. Then for every $J\subset I$ and $q\in\mathbb{N}$ one has $\dim_k(H_{dR}^q(]Y_J[^{\dagger}_{\cal X}))<\infty$, and also $\dim_k(H^q_{dR}(X^{\dagger}))<\infty$.
\end{satz}

{\sc Proof:} Since $X^{\dagger}=\cup_{i\in I}]Y_i[^{\dagger}_{\cal X}$ is a finite admissible covering, the second claim follows from the first. For the first we may, due to \ref{monovor}, assume that $J=I$, shrinking ${\cal X}$ if necessary. Passing to a finite covering, we may assume that ${\cal X}$ is affine, and that there is an \'{e}tale morphism ${\cal X}\stackrel{\phi}{\to}\spf(R<X_1,\ldots,X_n>/(X_1\ldots X_r-\pi))$ such that $Y_I=\cap_{i\in I}Y_i=\cap_{i=1}^rV(f_i)$ with $f_i=\phi^*(X_i)$. From the proof of \ref{monovor} we see that we may assume, again due to \ref{rigdrdf}, that there is a smooth $R$-algebra $\tilde{D}$, with weak completion $D^{\dagger}$, such that, if $B$ denotes the $k$-dagger algebra $D^{\dagger}\otimes_Rk$, we have an isomorphism of dagger spaces
$$]Y_I[^{\dagger}_{\cal X}=]V((\bar{f}_i)_{i=1,\ldots,r})[^{\dagger}_{\cal X}\cong$$$$\spm(B)\times\{x\in\spm(k<T_1,\ldots,T_r>^{\dagger}/(T_1\ldots T_r-\pi))|\quad|T_i(x)|<1\mbox{ for all }i\}$$such that $f_i$ corresponds to $T_i$. By \ref{semtub} this yields isomorphisms $$H_{dR}^q(]Y_I[^{\dagger}_{\cal X})\cong H_{dR}^q(\spm(B)\times V)\cong \oplus_{q=q_1+q_2}H_{dR}^{q_1}(\spm(B))\otimes_kH_{dR}^{q_2}(V)$$where we write $$V=\{(x_1,\ldots,x_{r-1})\in{\bf A}_k^{r-1}|\quad\prod_{i=1}^{r-1}|x_i|>|\pi|, |x_i|<1 \mbox{ for all }1\le i\le r-1\}.$$Now since $Y_I=\spec(\tilde{D}/(\pi))$ we have $H_{dR}^*(\spm(B))=H_{MW}^*(Y_I)$, and the latter is known to be finite dimensional by \cite{berfi},\cite{mebend}. But by \ref{semtub} also $H_{dR}^*(V)$ is finite dimensional. We are done.

\begin{lem}\label{faser} Let $n\ge s\ge r\ge 1$, let $$q:\spf(A)\to\spf(R<X_1,\ldots,X_n>/(X_1\ldots X_r-\pi))$$ be an \'{e}tale morphism and let $B=A/(q^*X_j)_{j=r+1,\ldots,s}$. There is an isomorphism$$]\spec(B/(\pi))[_{\spf(A)}\cong \spm(B\otimes_Rk)\times({\bf D}^0)^{s-r}$$such that the standard coordinates on $({\bf D}^0)^{s-r}$ correspond to $q^*X_{r+1},\ldots,q^*X_s$.
\end{lem}

{\sc Proof:} A strictly semi-stable formal $R$-scheme ${\cal X}$ carries a canonical log. structure: The log. structure ${\cal M}_{\cal X}$ associated with the divisor $X$, the reduction modulo $(\pi)$ of ${\cal X}$. In particular, $\spf(R)$ gives rise to a formal log. scheme $S$, and $({\cal X},{\cal M}_{\cal X})\to S$ is a log. smooth morphism of formal log. schemes. In our situation, $\spf(A)$ and $\spf(B)$ are strictly semi-stable formal $R$-schemes, and the embedding $({\cal U},{\cal M}_{\cal U})\to ({\cal V},{\cal M}_{\cal V})$ is an exact closed immersion of formal $S$-schemes.\\
Now we prove \ref{faser}, using the above log. structures. Because of \cite{EGA},IV, 18.3.2.1, we may work over the truncations mod $(\pi^n)$. Due to the extension property of morphisms from exact nilimmersions to log. smooth objects provided by \cite{kalo},3.11, as analogous to the classical extension property of morphisms from nilimmersions to smooth objects, one easily verifies the following transposition of \cite{EGA},0,19.5.4 to the log. context: Let $(\spec(B),M')\stackrel{g}{\to}(\spec(A),M)\stackrel{f}{\to}(\spec(C),N)$ be morphisms of affine log. schemes such that $g$ is an exact closed immersion defined by the ideal $I\subset A$ and such that $f$ and $f\circ g$ are log. smooth. Then $I/I^2$ is a projective $B$-module, and if $\hat{A}$ resp. $\hat{\cal S}$ are the respective $I$-adic completions of $A$ resp. $\sym_B(I/I^2)$, then there is a section $B\to\hat{A}$ together with an isomorphism of $B$-algebras $\hat{\cal S}\cong\hat{A}$. The lemma follows, because by the method \cite{berco},0.2.7 to compute tubes, all we have to do is to construct an $R$-isomorphism $B[[T_1,\ldots,T_{s-r}]]\cong\hat{A}$ with $T_i\mapsto q^*X_{i+r}$, where $\hat{A}$ is the $\ke(A\to B)$-adic completion of $A$.\\

\section{The finiteness theorem}

\addtocounter{satz}{1}{\bf \arabic{section}.\arabic{satz}}  We recall the terminology from \cite{dejo}. An $R$-variety is an integral separated flat $R$-scheme of finite type. Let $g:X\to\spec(R)$ be an $R$-variety and $X_i, i\in I$, the irreducible components of the closed fibre $X_{\bar {k}}$. For $J\subset I$ set $X_J=\cap_{i\in J}X_i$. 
Then $X$ is called strictly semi-stable over $R$ if the following (1)-(4) are fulfilled: (1) The generic fibre $X_k$ is smooth over $k$. (2) $X_{\bar {k}}$ is reduced, i.e. $X_{\bar {k}}=\cup_{i\in I}X_i$ scheme theoretically. (3) $X_i$ is a divisor on $X$ for all $i\in I$. (4) $X_J$ is smooth over $\bar{k}$ of codimension $|J|$ for all $J\subset I$, $J\ne\emptyset$.\\
Let $Z\subset X$ be Zariski closed with its reduced structure, such that $X_{\bar {k}}\subset Z$. Then $(X,Z)$ is called a strictly semi-stable pair over $R$ if the following (1)-(3) are fulfilled: (1) $X$ is strictly semi-stable over $R$. (2) $Z$ is a divisor with normal crossings on $X$. (3) Decompose $Z=Z_f\cup X_{\bar {k}}$ with $Z_f\to\spec(R)$ flat and let $Z_f=\cup_{i\in K}Z_i$ be the decomposition into irreducible components. Then $Z_L=\cap_{i\in L}Z_i$ is a union of strictly semi-stable $R$-varieties for all $L\subset K$, $L\ne\emptyset$.\\Note that for a strictly semi-stable pair $(X,Z)$ over $R$, the $\pi$-adic completion of $Z_f\to X$ is a strictly semi-stable formal pair over $R$ in the sense of \arabic{defsemist1}.\arabic{defsemist2}.\\

\addtocounter{satz}{1}{\bf \arabic{section}.\arabic{satz}} We call a dagger space $H$ quasi-algebraic, if there is an admissible covering of $H$ by dagger spaces $U$, which admit an open embedding into the dagger analytification of a $k$-scheme of finite type. If $H$ is quasi-algebraic, there is even an admissible covering of $H$ by open affinoids $U$ such that for each $U$ there are $n, r\in\mathbb{N}$, polynomials $f_j\in k[X_1,\ldots,X_n]$ and isomorphisms $U\cong\spm(k<X_1,\ldots,X_n>^{\dagger}/(f_1,\ldots,f_r))$ (see \cite{en1dag}, 2.18).

\begin{satz}\label{dejong}(\cite{dejo},6.5) If $Y$ is a proper $R$-variety and $Z\subset Y$ is a proper Zariski closed subset with $Y_{\bar {k}}\subset Z$, then there is a finite extension $R\to R'$ of complete discrete valuation rings, an $R'$-variety $X$, a proper surjective morphism of $R$-schemes $f:X\to Y$ and an open dense subscheme $U\subset Y$ such that $f^{-1}(U)\to U$ is finite and $(X,f^{-1}(Z)_{red})$ is a strictly semi-stable pair over $R'$.\\
\end{satz} 

\begin{satz} \label{endlich} For an affinoid quasi-algebraic dagger space $H$, the numbers $h_q^{dR}(H)$ are finite for all $q\ge0$.
\end{satz}

{\sc Proof:} Induction on $m=\dim(H)$. We may suppose $$H=\spm(k<X_1,\ldots,X_n>^{\dagger}/(f_1,\ldots,f_r))$$ with polynomials $f_j\in R[X_1,\ldots,X_n]$, and set $V=\spec(R[X_1,\ldots,X_n]/(f_1,\ldots,f_r))$. We regard $H$ as an open subspace of the dagger analytification of the generic fibre $V_k$ of $V$. The decomposition of $V$ into irreducible components induces a decomposition of $H$ into Zariski closed quasi-algebraic subspaces. Therefore we can reduce our claim by means of \arabic{eigensch1}.\arabic{eigensch2},(c) and the induction hypothesis to the case where $V$ is irreducible, $m=\dim(V_k)$. Since $h_q^{dR}(H)$ depends only on the reduced structure of $H$, we may assume $V$ is reduced, hence integral. 
Choose an open immersion $V\to Y$ into a projective $R$-variety $Y$. For the pair $(Y,Y_{\bar {k}})$ choose $R'$ and $f:X\to Y$ and $U\subset Y$ as in \ref{dejong} --- since $h_q^{dR}(H)$ is not effected by finite extensions of $k$, we may suppose $R=R'$. Let $X_k\stackrel{g_k}{\to}X'_k\stackrel{\pi_k}{\to}Y_k$ be the Stein factorization of the map $f_k:X_k\to Y_k$ of generic fibres; note that $X'_k$ is integral. From \cite{EGA},III, 4.4.1 it follows that there is a closed subscheme $T'\subset X'_k$, $\dim(T')<m$, such that for $T=g_k^{-1}(T')$ one has: $g_k|_{X_k-T}:(X_k-T)\to (X_k'-T')$ is an isomorphism. Since the locus of smoothness over $k$ is dense in $X'_k$ and in $Y_k$, and since $\pi_k$ is finite, we find a closed subscheme $S\subset Y_k$, $\dim(S)<m$, such that for $S'=X'_k\times_{Y_k}S$ one has: $Y_k-S$ and $X'_k-S'$ are smooth over $k$, and $(X'_k-S')\to(Y_k-S)$ is \'{e}tale. 
For a $k$-scheme $L$ of finite type, we denote by $L^{\dagger}$ its analytification as a dagger space. We regard $H$ as an open subspace of $Y_k^{\dagger}$ and set \\$H_X=H\times_{Y_k^{\dagger}}X_k^{\dagger},\quad\quad\quad\quad\quad H_{X'}=H\times_{Y_k^{\dagger}}X_k'^{\dagger},\quad\quad\quad H_S=H\times_{Y_k^{\dagger}}S^{\dagger},$\\$H_{X',S'}=H_{X'}\times_{X_k'^{\dagger}}S'^{\dagger},\quad\quad H_{X,T}=H_X\times_{X_k^{\dagger}}T^{\dagger},\quad\quad H_{X',T'}=H_{X'}\times_{X_k'^{\dagger}}T'^{\dagger}$.\\As auxiliary data we choose closed embeddings $H\to N$ resp. $H_{X'}\to M$ into smooth dagger spaces of pure dimension $n$ resp. $h$.\\We have $\dim_k(H^i(N,\dr(R{\underline\Gamma}_{*H}{\cal O}_N)))=h_{2n-i}^{dR}(H)$. By \ref{algana} there is a long exact sequence$$\ldots\to H^i(N,\dr(R{\underline\Gamma}_{*H_S}{\cal O}_N))\to H^i(N,\dr(R{\underline\Gamma}_{*H}{\cal O}_N))$$$$\to H^i(N-H_S,\dr(R{\underline\Gamma}_{*(H-H_S)}{\cal O}_{N-H_S}))\to H^{i+1}(\ldots.$$
The numbers $\dim_k(H^i(N,\dr(R{\underline\Gamma}_{*H_S}{\cal O}_N)))=h_{2n-i}^{dR}(H_S)$ are finite by induction hypothesis, so it is enough to show $\dim_k(H^i(N-H_S,\dr(R{\underline\Gamma}_{*(H-H_S)}{\cal O}_{N-H_S})))<\infty$, which by smoothness of $H-H_S$ is equivalent with $\dim_k(H^{i-2l}_{dR}(H-H_S))<\infty$ (where $l=\codim_N(H)$). The canonical maps $H^j_{dR}(H-H_S)\to H^j_{dR}(H_{X'}-H_{X',S'})$ are injective (\arabic{endspur1}.\arabic{endspur2}), so it is enough to show $\dim_k(H^j_{dR}(H_{X'}-H_{X',S'}))<\infty$. Since $H_{X'}-H_{X',S'}$ is smooth one has $$H^j_{dR}(H_{X'}-H_{X',S'})\cong H^{j+2t}(M-H_{X',S'},\dr(R{\underline\Gamma}_{*(H_{X'}-H_{X',S'})}{\cal O}_{M-H_{X',S'}}))$$ (where $t=\codim_M(H_{X'})$). By \ref{algana} there is a long exact sequence $$\ldots\to H^i(M,\dr(R{\underline\Gamma}_{*H_{X',S'}}{\cal O}_M))\to H^i(M,\dr(R{\underline\Gamma}_{*H_{X'}}{\cal O}_M))$$
$$\to H^i(M-H_{X',S'},\dr(R{\underline\Gamma}_{*(H_{X'}-H_{X',S'})}{\cal O}_{M-H_{X',S'}}))\to H^{i+1}(\ldots.$$Again $\dim_k(H^i(M,\dr(R{\underline\Gamma}_{*H_{X',S'}}{\cal O}_M)))=h_{2h-i}^{dR}(H_{X',S'})$, a finite number by induction hypothesis. It remains to show $\dim_k(H^i(M,\dr(R{\underline\Gamma}_{*H_{X'}}{\cal O}_M)))<\infty$.\\
Set $P=X_k^{\dagger}\times M$ and consider the closed immersion $H_X\to P$. From \ref{blasexa} we get (with $b=\dim(X_k)$) an exact sequence$$\ldots\to H^i(P,\dr(R{\underline\Gamma}_{*H_{X,T}}{\cal O}_P))\to H^{i-2b}(M,\dr(R{\underline\Gamma}_{*H_{X',T'}}{\cal O}_M))\oplus H^i(P,\dr(R{\underline\Gamma}_{*H_X}{\cal O}_P))$$$$\to H^{i-2b}(M,\dr(R{\underline\Gamma}_{*H_{X'}}{\cal O}_M))\to H^{i+1}(P,\dr(R{\underline\Gamma}_{*H_{X,T}}{\cal O}_P))\to\ldots.$$ Here $\dim_k(H^{2h-j}(M,\dr(R{\underline\Gamma}_{*H_{X',T'}}{\cal O}_M)))=h_j^{dR}(H_{X',T'})$ is finite by induction hypothesis. Furthermore $$H^{2h+i}(P,\dr(R{\underline\Gamma}_{*H_{X,T}}{\cal O}_P))=H^i(H_X,\dr(R{\underline\Gamma}_{*H_{X,T}}{\cal O}_{H_X}))$$since $H_X$ is smooth, and for all open affinoid $U\subset H_X$ also $$\dim_k(H^i(U,\dr(R{\underline\Gamma}_{*H_{X,T}}{\cal O}_{H_X})))=h_{2b-i}^{dR}(U\cap H_{X,T}),$$ a finite number by induction hypothesis. 
Since $H_X$ is quasi-compact it follows that $\dim_k(H^i(H_X,\dr(R{\underline\Gamma}_{*H_{X,T}}{\cal O}_{H_X})))<\infty$. So we are left with showing that $\dim_k(H^{2h+i}(P,\dr(R{\underline\Gamma}_{*H_X}{\cal O}_P)))=\dim_k(H_{dR}^i(H_X))$ is finite. By construction, the rigid space associated with $H$ is the generic fibre of an open formal subscheme of the $\pi$-adic completion of $Y$. It follows that the rigid space associated with $H_X$ is the generic fibre of an open formal subscheme of the $\pi$-adic completion of $X$. Thus $\dim_k(H_{dR}^i(H_X))<\infty$ follows from \ref{maxschni}. 

\begin{kor} \label{ende} If $X$ is a smooth quasi-compact dagger space and if $i:Z\to X$ is a closed immersion, then the $k$-vector spaces $H_{dR}^i(X-Z)$ are finite dimensional.
\end{kor}

{\sc Proof:} We may assume $X$ is affinoid. Choose a proper surjective morphism $f:\tilde{X}\to X$ such that $\tilde{X}$ is smooth, $\tilde{Z}=f^{-1}(Z)$ is a divisor with normal crossings on $\tilde{X}$ and $(\tilde{X}-\tilde{Z})\to(X-Z)$ is an isomorphism (\ref{hironaka}). Passing to an appropriate finite affinoid admissible open covering of $\tilde{X}$ we see that we may assume from the beginning: $X$ is affinoid and $Z$ is a normal crossings divisor on $X$ such that all its irreducible components are smooth. Now note that if $W=X$ or if $W$ is the intersection of some irreducible components of $Z$, then $W$ is quasi-algebraic: Indeed, since it it smooth and affinoid, it follows from \cite{elk}, th. 7, p. 582 that the associated rigid space $W'$ can be defined by polynomials,. In particular $W'$ is the rigid space associated to a quasi-algebraic affinoid dagger space $W_1$, and by \cite{en1dag}, 1.15 there exists a (non-canonical) isomorphism $W\cong W_1$. Now \ref{endlich} says that $H^q(X,\dr(R{\underline\Gamma}_{*W}{\cal O}_{X}))$ is finite dimensional for all $q$ and all such $W$. Repeated application of \ref{algana} gives the Corollary.\\ 
 
\begin{satz}\label{kokomend} Let $X$ be a quasi-compact smooth dagger space, $U\subset X$ a quasi-compact admissible open subset. Then $H_{dR}^q(X-U)$ is finite dimensional for all $q\in\mathbb{N}$.
\end{satz}

{\sc Proof:} We prove by induction on $n\in\mathbb{N}$:\\${\bf (a_n)}$ For every quasi-compact smooth dagger space $X$ with $\dim(X)\le n$, every quasi-compact open $U\subset X$ and every $q\in\mathbb{N}$, we have $\dim_k(H^q_{dR}(X-U))<\infty$.\\
${\bf (b_n)}$ For every quasi-compact smooth dagger space $Y$, every closed immersion $T\hookrightarrow Y$ with $\dim(T)\le n$, every quasi-compact open $V\subset Y$ and every $q\in\mathbb{N}$, we have $\dim_k(H^q(Y-V,\dr(R{\underline\Gamma}_{*T}{\cal O}_Y)))<\infty.$\\
Here ${\bf (a_0)}$ is evident, and so is ${\bf (b_0)}$ because of \ref{dmodgys}.\\
${\bf (b_{n-1})\Rightarrow(a_n):}$ We may suppose $\dim(X)=n$ and $X$ is affinoid and connected. Since $U$ is the union of finitely many rational subdomains of $X$, Mayer-Vietoris sequences allow us to reduce to the case where $U$ is a rational subdomain of $X$. Choose an affine formal $R$-scheme ${\cal X}$ such that its generic fibre ${\cal X}_k$ is the rigid space associated with $X$. For subschemes $Z\subset{\cal X}_{\bar {k}}$ we denote by $]Z[_{\cal X}^{\dagger}$ the open dagger subspace of $X$ corresponding to $]Z[_{\cal X}\subset {\cal X}_k$. By \cite{bolu}, after an admissible blowing up and further localization we may suppose: There is a closed subscheme $Z\subset{\cal X}_{\bar {k}}$ such that $]Z[_{\cal X}^{\dagger}=X-U$. From \cite{elk}, th. 7, p. 582 it follows that ${\cal X}$ is locally defined by polynomials, i.e. we may suppose $${\cal X}=\spf(R<X_1,\ldots,X_m>/(f_1,\ldots,f_r))$$ with $f_i\in R[X_1,\ldots,X_m]$. We view $X$ as an open subspace of the dagger analytification of the generic fibre $T_k$ of $T=\spec(R[X_1,\ldots,X_m]/(f_1,\ldots,f_r))$. Since $X$ is smooth, it is contained in the dagger analytification of a single irreducible component of $T_k$. This component is the generic fibre of a closed subscheme of $T$; dividing out the $\pi$-torsion and the nilpotent elements of its coordinate ring, we see that we may suppose that $T$ is integral. We may also suppose that $Z$ is defined by a single equation in $T_{\bar{k}}={\cal X}_{\bar {k}}$. We choose a closed subscheme $Y\subset T$ defined by a single equation such that $Y_{\bar {k}}=Z$. Furthermore we choose an open embedding $T\to\bar{X}$ into a projective $R$-variety $\bar{X}$, and define $\bar{Y}$ to be the schematic closure of $Y$ in $\bar{X}$. Since ${\bf (a_n)}$ is proven in case $U=\emptyset$ by \ref{endlich}, we may suppose $\bar{Y}\cup\bar{X}_{\bar {k}}$ is a proper subset of $\bar{X}$. Therefore we can apply \ref{dejong} to $(\bar{X},\bar{Y}\cup\bar{X}_{\bar {k}})$: Performing a base change with a finite extension of $R$, we may assume that there is a surjective proper morphism $\bar{\phi}:\bar{V}\to\bar{X}$ of $R$-varieties, an open dense subscheme of $\bar{X}$ over which $\bar{\phi}$ is finite, and such that $(\bar{V},\bar{\phi}^{-1}(\bar{Y}\cup\bar{X}_{\bar{k}})_{red})$ is a strictly semi-stable pair over $R$. 
Let $\bar{\cal V}$, resp. $\bar{\cal W}$, resp. $\bar{\cal X}$ be the $\pi$-adic formal completion of $\bar{V}$, resp. $\bar{\phi}^{-1}(\bar{Y})_{red}$, resp. $\bar{X}$, and set ${\cal V}=\bar{\cal V}\times_{\bar{\cal X}}{\cal X}$ and ${\cal W}=\bar{\cal W}\times_{\bar{\cal V}}{\cal V}$. Let $\bar{V}^{\dagger}_k$ be the dagger analytification of the generic fibre of $\bar{V}$, and let $V$ be its open dagger subspace whose associated rigid space is ${\cal V}_k$. As before, for subschemes $S\subset{\cal V}_{\bar{k}}$, we denote by $]S[^{\dagger}_{\cal V}$ the open dagger subspace of $V$ corresponding to $]S[_{\cal V}\subset {\cal V}_k$. Now $\bar{\phi}$ induces a morphism $\phi_k:V\to X$ of dagger spaces, and if $F=X-U$, we have $]{\cal W}_{\bar {k}}[^{\dagger}_{\cal V}=\phi_k^{-1}(F)$.\\
Claim: It is enough to show $\dim_k(H_{dR}^*(]{\cal W}_{\bar {k}}[^{\dagger}_{\cal V}))<\infty$ for all $q\in\mathbb{N}$.\\
Using the induction hypothesis ${\bf (b_{n-1})}$, this can be shown similarly as in the proof of \ref{endlich}: Let $V\stackrel{f}{\to}D\stackrel{g}{\to}X$ be the Stein factorization of $\phi_k$ (for example obtained from the Stein factorization of the algebraic morphism $\bar{\phi}$). Let $Q=g^{-1}(F)$, let $D\to P'$ be a closed immersion into a smooth affinoid dagger space $P'$ and let $P\subset P'$ be an open subspace such that $P'-P$ is quasi-compact and open in $P'$ and such that $Q\to P'$ factorizes through a closed immersion $Q\to P$. Let  $T\to F$ be a closed immersion with $\dim(T)<n$, such that for $L=g^{-1}(T)$ we have: $Q-L$ is smooth and $(Q-L)\to(F-T)$ is \'{e}tale. By \ref{algana} there are long exact sequences
$$\ldots\to H^i(F,\dr(R{\underline\Gamma}_{*T}{\cal O}_F))\to H_{dR}^i(F)$$$$\to H_{dR}^i(F-T)\to H^{i+1}(F,\dr(R{\underline\Gamma}_{*T}{\cal O}_F))\to\ldots$$
and$$\ldots\to H^i(P,\dr(R{\underline\Gamma}_{*L}{\cal O}_P))\to H^i(P,\dr(R{\underline\Gamma}_{*Q}{\cal O}_P))$$$$\to H^i(P-L,\dr(R{\underline\Gamma}_{*Q-L}{\cal O}_{P-L}))\to H^{i+1}(P,\dr(R{\underline\Gamma}_{*L}{\cal O}_P))\to\ldots.$$Now
$H^i(F,\dr(R{\underline\Gamma}_{*T}{\cal O}_F))$ and $H^i(P,\dr(R{\underline\Gamma}_{*L}{\cal O}_P))$ are finite dimensional by induction hypothesis ${\bf (b_{n-1})}$. Since $Q-L$ is smooth we have $$H^i(P-L,\dr(R{\underline\Gamma}_{*Q-L}{\cal O}_{P-L}))\cong H_{dR}^{i-2l}(Q-L)$$ (where $l=\codim_P(Q)$), and by \arabic{endspur1}.\arabic{endspur2} the canonical maps $H^q_{dR}(F-T)\to H_{dR}^q(Q-L)$ are injective. 
Together we obtain that to prove $\dim_k(H^q_{dR}(F))<\infty$, it is enough to prove $\dim_k(H^q(P,\dr(R{\underline\Gamma}_{*Q}{\cal O}_P)))<\infty$ for all $q\in\mathbb{N}$. Now choose a closed immersion $G\to Q$ with $\dim(G)<n$ such that for $H=G\times_Q]{\cal W}_{\bar {k}}[^{\dagger}_{\cal V}$ we have: $(]{\cal W}_{\bar {k}}[^{\dagger}_{\cal V}-H)\to(Q-G)$ is an isomorphism. 
Let $E=\bar{V}_k^{\dagger}\times P$ with $]{\cal W}_{\bar {k}}[^{\dagger}_{\cal V}$ diagonally embedded. By \ref{blasexa} there is a long exact sequence$$\ldots\to H^i(E,\dr(R{\underline\Gamma}_{*H}{\cal O}_E))\to H^{i-2n}(P,\dr(R{\underline\Gamma}_{*G}{\cal O}_P))\oplus H^i(E,\dr(R{\underline\Gamma}_{*]{\cal W}_{\bar {k}}[^{\dagger}_{\cal V}}{\cal O}_E))$$
$$\to H^{i-2n}(P,\dr(R{\underline\Gamma}_{*Q}{\cal O}_P))\to H^{i+1}(E,\dr(R{\underline\Gamma}_{*H}{\cal O}_E))\to\ldots.$$We have $\dim_k(H^i(E,\dr(R{\underline\Gamma}_{*H}{\cal O}_E)))<\infty$ and $\dim_k(H^j(P,\dr(R{\underline\Gamma}_{*G}{\cal O}_P)))<\infty$ by induction hypothesis ${\bf (b_{n-1})}$. On the other hand $$H^i(E,\dr(R{\underline\Gamma}_{*]{\cal W}_{\bar {k}}[^{\dagger}_{\cal V}}{\cal O}_E))=H_{dR}^{i-2n}(]{\cal W}_{\bar {k}}[^{\dagger}_{\cal V}),$$ and altogether the claim follows.\\
Now we prove $\dim_k(H_{dR}^q(]{\cal W}_{\bar {k}}[^{\dagger}_{\cal V}))<\infty$. Note that ${\cal W}$ is the $\pi$-adic formal completion of $\bar{\phi}^{-1}(Y)_{red}$. Let ${\cal W}={\cal W}_f\cup{\cal W}_0$ be the decomposition into the $R$-flat part ${\cal W}_f$ and the $\pi$-torsion part ${\cal W}_0$. Since $Y$ is defined by a single equation in $T$, the same is true for $\bar{\phi}^{-1}(Y)_{red}$ in $\bar{\phi}^{-1}(T)_{red}$. Since $(\bar{V},\bar{\phi}^{-1}(\bar{Y}\cup\bar{X}_{\bar {k}})_{red})$ is a strictly semi-stable pair over $R$, this means that $\bar{\phi}^{-1}(Y)_{red}$ is the union of some irreducible components of $\bar{\phi}^{-1}(Y\cup T_{\bar {k}})_{red}$. After passing to a finite Zariski open covering of ${\cal V}$ we may therefore suppose: There is an \'{e}tale morphism $${\cal V}\stackrel{\psi}{\to}\spf(R<X_1,\ldots,X_m>/(X_1\ldots X_r-\pi))$$ for some $m\ge r\ge1$, and if we set ${\cal W}^i=V(\psi^*X_i)$ for $i\le m$, there are subsets $J\subset\{r+1,\ldots,m\}$ and $N\subset\{1,\ldots,r\}$ such that ${\cal W}_f=\cup_{i\in J}{\cal W}^i$ and ${\cal W}_0=\cup_{i\in N}{\cal W}^i$. The covering $]{\cal W}_{\bar {k}}[^{\dagger}_{\cal V}=\cup_{i\in J\cup N}]{\cal W}^i_{\bar {k}}[^{\dagger}_{\cal V}$ is admissible, so it is enough to show: 
For all $I\subset J$, all $M\subset N$, all $q\in\mathbb{N}$ we have $$\dim_k(H_{dR}^q(\cap_{i\in I\cup M}]{\cal W}^i_{\bar {k}}[^{\dagger}_{\cal V}))<\infty.$$ Set ${\cal R}=\cap_{i\in I}{\cal W}^i$, a strictly semi-stable formal $R$-scheme. By construction, the equations defining the closed immersion ${\cal R}\to{\cal V}$ are contained in ${\cal O}_V$; they define a Zariski closed dagger subspace of $V$ whose associated rigid space is the generic fibre of ${\cal R}$. Similar as before we define its open subspaces $]S[^{\dagger}_{\cal R}$ for $S\subset {\cal R}_{\bar {k}}$. Let $C={\cal R}_{\bar {k}}\cap(\cap_{i\in M}{\cal W}^i)$. We claim $$H_{dR}^q(\cap_{i\in I\cup M}]{\cal W}^i_{\bar {k}}[^{\dagger}_{\cal V})=H_{dR}^q(]C[^{\dagger}_{\cal V})\cong H_{dR}^q(]C[^{\dagger}_{\cal R}).$$ Indeed, from \ref{faser} we derive an isomorphism $]C[_{\cal R}\times({\bf D}^0)^{|I|}\cong]C[_{\cal V}$ of rigid spaces. Because of \ref{rigdrdf}, to prove our claim we may therefore assume that there is an isomorphism of dagger spaces$$]C[^{\dagger}_{\cal R}\times({\bf D}^0)^{|I|}\cong]C[^{\dagger}_{\cal V}$$ and then the claim is obvious. But $H_{dR}^q(]C[^{\dagger}_{\cal R})$ is finite dimensional by \ref{maxschni}, because $C$ is the intersection of some irreducible components of ${\cal R}_{\bar {k}}$.\\
${\bf (a_n)+(b_{n-1})\Rightarrow(b_n):}$ Passing to a finite affinoid covering of $Y$ we may suppose (\ref{hironaka}): There is a smooth $P$ and a proper surjective $P\stackrel{f}{\to}Y$ such that $S=f^{-1}(T)$ is a divisor with normal crossings on $P$ and $(P-S)\to(Y-T)$ is an isomorphism. 
As in ${\bf (b_{n-1})\Rightarrow(a_n)}$ one shows that it suffices to prove $\dim_k(H^q(f^{-1}(Y-V),\dr(R{\underline\Gamma}_{*S}{\cal O}_P)))<\infty$; so we suppose $T$ is a divisor with normal crossings on $Y$. But then, in view of \arabic{eigensch1}.\arabic{eigensch2}(c), the problem is equivalent with the one where $T$ is smooth (of arbitrary codimension). Thus it is reduced by means of \arabic{eigensch1}.\arabic{eigensch2}, \ref{dmodgys} to the induction hypothesis ${\bf (a_n)}$.\\

\addtocounter{satz}{1}{\bf \arabic{section}.\arabic{satz}} As an application of \ref{kokomend}, we will show in \cite{en3dag} that the de Rham cohomology groups of smooth rigid Stein spaces are topologically separated for their canonical topology, hence are Fr\'{e}chet spaces. We also define reasonable de Rham cohomology groups for arbitrary rigid spaces 
(the underlying idea is that of \ref{rigdrdf}), and derive from \ref{kokomend} K\"unneth and duality formulas for them.  

\begin{kor}\label{rifi} For a $\bar{k}$-scheme $Y$ of finite type, the $k$-vector spaces $H_{rig}^q(Y/k)$ are finite dimensional for all $q\in\mathbb{N}$. 
\end{kor}

{\sc Proof:} Let ${\cal X}$ be a proper smooth $\pi$-adic tf formal $\spf(R)$-scheme, $Y\to{\cal X}_{\bar {k}}$ an immersion with schematic closure $j:Y\to\bar{Y}$ in ${\cal X}_{\bar {k}}$. Then $]\bar{Y}[_{\cal X}$ is a partially proper rigid space, therefore equivalent with a dagger space $Q$. 
Let $X$ be the open subspace of $Q$ whose underlying set is identified with $]Y[_{\cal X}$. From \cite{en1dag}, 5.1 we get an isomorphism $H_{rig}^q(Y/k)\cong H^q_{dR}(X)$, but $H^q_{dR}(X)$ is finite dimensional by \ref{kokomend}. 


\end{document}